\documentclass[11pt, twoside]{amsart}

\usepackage[utf8]{inputenc}
\usepackage[T1]{fontenc}
\usepackage{amssymb,amsmath,amstext}
\usepackage{hyperref}

\newtheorem{theor}{Theorem}[section]
\newtheorem{lem}[theor]{Lemma}
\newtheorem{defin}[theor]{Definition}

\newtheorem{prop}[theor]{Proposition}

\newtheorem{cor}[theor]{Corollary}
\newtheorem{rem}[theor]{Remark}

\newtheorem{fact}[theor]{Fact}

\newtheorem{observation}[theor]{Observation}

\numberwithin{equation}{section}

\newcommand{\acl}{\mathrm{acl}}
\newcommand{\dcl}{\mathrm{dcl}}
\newcommand{\tp}{\mathrm{tp}}

\newcommand{\es}{\emptyset}

\newcommand{\nts}{\negthickspace}
\newcommand{\uhrc}{\nts \upharpoonright \nts}

\newcommand{\meq}{^{\mathrm{eq}}}

\newcommand{\mcA}{\mathcal{A}}
\newcommand{\mcB}{\mathcal{B}}
\newcommand{\mcC}{\mathcal{C}}
\newcommand{\mcD}{\mathcal{D}}
\newcommand{\mcE}{\mathcal{E}}
\newcommand{\mcF}{\mathcal{F}}
\newcommand{\mcG}{\mathcal{G}}
\newcommand{\mcH}{\mathcal{H}}

\newcommand{\mcK}{\mathcal{K}}

\newcommand{\mcM}{\mathcal{M}}
\newcommand{\mcN}{\mathcal{N}}

\newcommand{\mcP}{\mathcal{P}}

\newcommand{\mbC}{\mathbf{C}}
\newcommand{\mbD}{\mathbf{D}}

\newcommand{\mbF}{\mathbf{F}}

\newcommand{\mbK}{\mathbf{K}}
\newcommand{\mbS}{\mathbf{S}}
\newcommand{\mbT}{\mathbf{T}}

\newcommand{\mrqf}{\mathrm{qf}}

\newcommand{\ind}{\raisebox{-2pt}[5pt][0pt]{$\smile$} \hspace*{-6.8pt}\raisebox{3pt}[5pt][0pt]{$|$} \; \: }
\newcommand{\nind}{\raisebox{-2pt}[5pt][0pt]{$\smile$} 
\hspace*{-6.8pt}\raisebox{3pt}[5pt][0pt]{$|$}\hspace*{-6.8pt}
\raisebox{3pt}[5pt][0pt]{$\diagup$} }

\newcommand{\rng}{\mathrm{rng}}

\title[On constraints and dividing]
{On constraints and dividing \\ in ternary homogeneous structures}

\author{Vera Koponen}

\address{Vera Koponen, Department of Mathematics, Uppsala University, Box 480,
75106 Uppsala, Sweden.}

\email{vera.koponen@math.uu.se}

\date{17 July, 2018. Revised version}

\begin{document}

\begin{abstract}
Let $\mcM$ be ternary, homogeneous and simple. 
We prove that if $\mcM$ is finitely constrained, then it is supersimple with finite SU-rank 
and dependence is $k$-trivial for some $k < \omega$ and for finite sets of real elements. 
Now suppose that, in addition, $\mcM$ is supersimple with SU-rank~1. 
If $\mcM$ is finitely constrained then algebraic closure in $\mcM$ is trivial.
We also find connections between the nature of the constraints of $\mcM$,
the nature of the amalgamations allowed by the age of $\mcM$, 
and the nature of definable equivalence relations.
A key method of proof is to ``extract'' constraints (of $\mcM$) from
instances of dividing and from definable equivalence relations.
Finally, we give new examples, including an uncountable family, 
of ternary homogeneous supersimple structures of SU-rank~1.

\noindent
{\em Keywords}: model theory, homogeneous structure, simple theory, constraint, dividing, amalgamation.
\end{abstract}

\maketitle

\section{Introduction}

\noindent
We call a structure $\mcM$  {\em homogeneous} if it is countable, has a finite relational vocabulary and every isomorphism between 
finite substructures of $\mcM$ can be extended to an automorphism of $\mcM$.
Homogeneous structures are of interest in various areas, including Ramsey theory, constraint satisfaction problems, permutation group theory 
and topological dynamics; surveys include~\cite{BP, HN, Mac11, Nes}.
From the point of view of pure model theory they are interesting because they can also be characterized as
the countable structures (with finite relational vocabulary) which have elimination of quantifiers, or
as the countable structures (with finite relational vocabulary) which are so-called Fra\"{i}ss\'{e} limits of 
``amalgamation classes'' of finite structures (see for example \cite[Chapter~7]{Hod}).
The later two characterizations give homogeneous structures a rather ``concrete'' character.
Nevertheless, there are uncountably many homogeneous digraphs~\cite{Hen72}.
Some classes of homogeneous structures with additional properties, such as for example 
homogeneous (di)graphs, partial orders and stable (infinite) structures
have been classified \cite{Che98, LW, Schm, Lach97}.

It is of course tempting to try to classify other classes of homogeneous structures.
A natural direction is to consider simple homogeneous structures, or some other class of homogeneous structures 
in which there is a reasonably well behaved notion of independence.
What makes homogeneous structures intriguing to me is that, although they are fairly concrete
(and useful for providing examples),
typical model theoretic questions such as whether a simple homogeneous structure is supersimple, 
has finite SU-rank, or has a nontrivial pregeometry, turn out to be challenging.
Moreover, the class of homogeneous simple structures can be seen as the most uncomplicated kind of $\omega$-categorical simple structures
in the sense that in a homogeneous structure every definable relation is definable by
a quantifier-free formula (a ``$\Delta_0$-formula''), while in an $\omega$-categorical structure a definable relation may be
definable only by a formula of higher complexity, in terms of, say, quantifier alternations.
For example, there is an $\omega$-categorical supersimple structure with SU-rank~1 in which algebraic closure induces
a non-locally modular pregeometry, as shown first by Hrushovski~\cite{Hru03, Kim_book}.
But we do not know yet if a similar example exists if `$\omega$-categorical' is replaced with `homogeneous'. 
My guess is that the answer is `no', and if this is correct one can ask, given any $k < \omega$,
if there is an $\omega$-categorical supersimple structure with SU-rank~1 in which algebraic closure induces
a non-locally modular pregeometry and all definable relations are defined by $\Sigma_k$ (or $\Pi_k$) formulas.
In general, for a given phenomenon, one can ask what the minimal ``definitional complexity'' of relations must be in a structure
for that phenomenon to appear in it.

A relational structure (i.e. one with relational vocabulary) will be called {\em $k$-ary} if no relation symbol of its vocabulary has arity higher than $k$. 
We say {\em binary} and {\em ternary} instead of 2-ary and 3-ary, respectively.
The fine structure of all binary simple homogeneous structures is well understod in terms of supersimplicity, SU-rank,
pregeometries, the behaviour of dividing and
the nature of definable sets of SU-rank~1 \cite{AK, Kop16PAM, KopPrimitive, KopSub}. 
However, the arguments for binary structures do not carry over to ternary structures. 
This is at least partly to be expected, because at least one new phenomenon appear when passing from binary to ternary:
Every binary supersimple homogeneous structure with SU-rank~1 is a random structure
(in the sense of Definition~\ref{definition of age, permitted etc} below), 
but there are plenty of ternary supersimple structures with SU-rank~1 which are not random, 
as witnessed by the examples in Section~\ref{Examples}.

The notion of $n$-complete amalgamation in the sense of \cite[Definition~2.2]{PalRecent} is related to 
random structures (in the presence of homogeneity). The independence theorem of simple theories is equivalent to 3-complete amalgamation.
In the theory of binary simple homogeneous structures the combination of the independence theorem and the binarity (of the structures)
plays a crucial role. Palac\'{i}n \cite{PalRecent} has recently shown that if a $k$-ary relational structure
is homogeneous, supersimple and has $(k+1)$-complete amalgamation, then it is a random structure.
Kruckman~\cite[sections~5.3 and~5.5]{KruPhD} has also investigated ``higher dimensional'' amalgamation properties 
in the context of $\omega$-categorical (possibly simple) structures.
Consequently, we will be concerned with ternary homogeneous simple structures which do not necessarily
have 4-complete amalgamation.
Conant has recently studied the relationship of the free amalgamation property (of a class of finite structures) to
dividing, rosiness, simplicity and related notions (of the Fra\"{i}ss\'{e} limit). 
Both Conant \cite{Con} and Palac\'{i}n \cite{PalRecent} have proved results which relate
the free amalgamation property to (super)simplicity and SU-rank one.
The relevance of definable (with parameters) equivalence relations in a context related to the one of this article has
been realized earlier by Kruckman~\cite[Section~5.5]{KruPhD} (we will deal with equivalence relations in
Section~\ref{Definable equivalence relations and weakly isolated contraints}). 
Akhtar and Lachlan \cite{AL} have studied homogeneous (not necessarily simple) 3-hypergraphs. 
They have proved results about the age of such a 3-hypergraph, and, in particular, characterize the constraint if there is only one.

Needless to say, if the study of simple homogeneous structures requires a special treatment of 
$k$-ary simple homogeneous structures for every $k < \omega$,
then the study will never be completed.
However, I hope that for some small $k$ (perhaps $k = 4$ or $k=5$) the arguments for the case $k$ work out
for all $k' \geq k$, at least with respect to general questions such as whether a structure is supersimple, has finite SU-rank, 
or what the nature of definable pregeometries is.
Some justification for my optimism comes from the theory of smoothly approximable structures~\cite{CH}
and from Corollaries 5.3 and~5.4 in~\cite{KopSub} where the nature of $4$-types over $\es$ is crucial.

In this article we will study the interplay between, on the one hand, constraints of a ternary homogeneous structure
and, on the other hand, the behaviour of dividing (and indirectly the algebraic closure) and the existence 
of nontrivial definable (with parameters) equivalence relations on the set of realizations of 
a complete nonalgebraic 1-type (over the same parameters).
The constraints of a homogeneous structure are the ``minimal'' structures (for the same vocabulary) which can not
be embedded into it (see Definition~\ref{definition of age, permitted etc}).
The first main result (Theorem~\ref{the rank is finite})
tells that if a structure is ternary, homogeneous, finitely constrained and simple, 
then it is supersimple with finite SU-rank and dependence is, for some $k < \omega$, $k$-trivial for finite sets of real elements.
As a background recall that all stable homogenous  structures are finitely constrained~\cite[Theorem~5]{Lach86}, but as we will see in 
Section~\ref{Uncountably many ternary homogeneous simple structures}, 
there are ternary homogeneous simple structures which are not finitely constrained.
Then we see that if $\mcM$ is ternary, homogeneous, finitely constrained and supersimple with SU-rank~1, 
then algebraic closure and dependence in $\mcM$ are trivial (Theorem~\ref{trivial acl if the rank is one}).

In Section~\ref{Definable equivalence relations and weakly isolated contraints}
we turn to the special case of SU-rank~1, that is, we study ternary homogeneous (not necessarily finitely constrained)
supersimple structures of SU-rank~1.
Here we show that the existence of nontrivial definable
equivalence relations on the set of realizations of a nonalgebraic 1-type over a finite set
plays a crucial role for undestanding the nature of the constraints of the structure and the nature of its algebraic closure,
and vice versa.
We also see that if the age of the structure has the free amalgamation property, then there are no definable
equivalence relations on any nonalgebraic types and all constraints are of a particular kind.

Section~\ref{Definability of binary random structures}
elaborates more on the topic of nontrivial definable equivalence relations and shows that, under
some extra conditions, no nontrivial equivalence relation is definable with only one parameter.
It follows, under the extra conditions, that a binary random structure is definable (with only one parameter)
in the structure that we started with. If the extra conditions are not satisfied, then we only conclude that
a binary structure $\mcN$ is definable in the original structure $\mcM$ such that 
$\mcN$ is random relative to an equivalence relation on its universe which is $\es$-definable in $\mcN$.

Section~\ref{Examples} gives examples of ternary homogeneous supersimple structures with SU-rank~1 and
degenerate algebraic closure.
The examples in sections~\ref{Uncountably many ternary homogeneous simple structures} 
and~\ref{A ternary homogeneous structure which is definable in the generic tournament}
seem to be new.
Section~\ref{Uncountably many ternary homogeneous simple structures}
shows that, with a ternary relation symbol $R$ and vocabulary $V = \{R\}$, there are
uncountably many nonisomorphic $V$-structures that are homogeneous and supersimple with SU-rank~1
and degenerate algebraic closure.
The last section gives a list of problems.

\section{Preliminaries}

\subsection{Notation and general concepts}

Structures are denoted by calligraphic letters $\mcA, \mcB, \ldots, \mcM, \mcN$ and their universes by
the corresponding non-calligraphic letters $A, B, \ldots, M, N$.
Usually (but not always) infinite structures are denoted by $\mcM$ or $\mcN$, possibly with indices, and finite structures
by $\mcA, \mcB, \ldots,$ possibly with indices.
The complete theory of a structure $\mcM$ is denoted by $Th(\mcM)$.
If $V$ is a vocabulary, $\mcM$ a $V$-structure and $W \subseteq V$, then
$\mcM \uhrc W$ denotes the reduct of $\mcM$ to $W$.
If $V$ is a relational vocabulary, $\mcM$ is a $V$-structure and $A \subseteq M$, then $\mcM \uhrc A$ denotes
the substructure of $\mcM$ with universe $A$.
By $\bar{a}, \bar{b}, \ldots$ we denote finite sequences/tuples of elements.
If $A$ is a set then `$\bar{a} \in A$' usually means that $\bar{a}$ is a finite sequence of elements from $A$.
The length of $\bar{a}$ is denoted by $|\bar{a}|$ and we may write $\bar{a} \in A^n$ if we want to emphasize that
$\bar{a}$ is a sequence of length $n$ all of which elements belong to $A$.
By $\rng(\bar{a})$ we denote the set of elements that occur in $\bar{a}$.
For sequences $\bar{a}$ and $\bar{b}$, $\bar{a}\bar{b}$ denotes the concatenation of them.
Sometimes we abuse notation and write `$\bar{a}$' in instead of `$\rng(\bar{a})$' and `$AB$' instead of `$A \cup B$'.
For a formula $\varphi(\bar{x})$, or type $p(\bar{x})$, $\varphi(\mcM)$, respectively $p(\mcM)$, denotes
the set of tuples of elements from $\mcM$ which satisfy/realize it.

Given a structure $\mcM$ and $A \subseteq M$, $S_n^\mcM(A)$ denotes the set of complete $n$-types over $A$
(with respect to $\mcM$).
For a structure $\mcM$ and $\bar{a} \in M$, $\tp_\mcM^\mrqf(\bar{a})$ denotes the set of  quantifier
free formulas satisfied by $\bar{a}$ in $\mcM$
(while $\tp_\mcM(\bar{a})$, as usual, denotes the complete type of $\bar{a}$ in $\mcM$ and
$\tp_\mcM(\bar{a} / B)$ denotes the complete type of $\bar{a}$ over $B$ in $\mcM$ if $B \subseteq M$).
We sometimes write $\bar{a} \equiv_\mcM \bar{b}$, or $\bar{a} \equiv_\mcM^\mrqf \bar{b}$ instead of
$\tp_\mcM(\bar{a}) = \tp_\mcM(\bar{a})$, or $\tp_\mcM^\mrqf(\bar{a}) = \tp_\mcM^\mrqf(\bar{b})$, respectively.

We assume familiarity with the notions of dividing and forking, as well as simple theories and SU-rank,
as can be found in \cite{Cas, Kim_book}, for example.
When saying that a structure $\mcM$ is simple, supersimple, or has finite SU-rank, then we mean that its complete
theory, denoted $Th(\mcM)$, has the corresponding property.
In some arguments we will consider dividing in different structures and in this context we may write `$\ind^\mcM$' to
indicate that we consider dividing in $\mcM$.
If we say that a structure is simple, then we assume that it is infinite.

A structure $\mcM$ is said to be {\em $k$-transitive} if $Th(\mcM)$ has only one complete $k$-type over $\es$
which implies that all $k$ elements are different.
An equivalence relation is called {\em nontrivial} if it has at least two equivalence classes and at least one equivalence
class contains more than one element.
Suppose that $\mcM$ and $\mcN$ are structures with the same universe but possibly (and typically) with different vocabularies.
We say that $\mcN$ is a {\em reduct} of $\mcM$ if, for every $n < \omega$ and $R \subseteq N^n$,
if $R$ is $\es$-definable in $\mcN$ then it is $\es$-definable in $\mcM$.

\subsection{Classes of finite structures}\label{Classes of finite structures}

We assume familiarity with the basic theory of ``amalgamation classes'' of finite structures and Fra\"{i}ss\'{e} limits
(as explained in \cite{Hod} for example), but nevertheless some terminology is defined below to avoid confusion.

\begin{defin}\label{amalgamation properties}{\rm
Let $V$ be a finite relational vocabulary.\\
(i) Let $\mcM$ be a $V$-structure and $\bar{a} \in M$. For $R \in V$ we say that $\bar{a}$ is an {\em $R$-relationship (in $\mcM$)}
if $\bar{a} \in R^\mcM$. We say that $\bar{a}$ is a {\em relationship (in $\mcM$)} if for some $R \in V$, 
$\bar{a}$ is an $R$-relationship.\\
(ii) If $\mcA$ and $\mcB$ are two $V$-structures, then the {\em free amalgam} of $\mcA$ and $\mcB$ is the unique
structure $\mcC$ such that $C = A \cup B$ and for every $R \in V$ and every tuple $\bar{d} \in C$,
$\bar{d} \in R^\mcC$ if and only if $\bar{d} \in R^\mcA \cup R^\mcB$. 
If $\mcE$ is a substructure of $\mcA$ and of $\mcB$ of maximal cardinality, then we may also say that $\mcC$
(as defined above) is the {\em free amalgam of $\mcA$ and $\mcB$ over $\mcE$}. \\
(iii) A class $\mbK$ of finite $V$-structures has the {\em hereditary property} if it is closed under substructures
(i.e. if $\mcA \subseteq \mcB \in \mbK$, then $\mcA \in \mbK$). \\
(iv) Suppose that $\mbK$ is a class of finite $V$-structures which is closed under isomorphism.
Then $\mbK$ has the
\begin{itemize}
\item[(a)] {\em amalgamation property} if whenever $\mcA, \mcB, \mcC \in \mbK$ and 
$f_\mcB: \mcA \to \mcB$ and $f_\mcC : \mcA \to \mcC$ are embeddings, then there are $\mcD \in \mbK$ and embeddings 
$g_\mcB : \mcB \to \mcD$ and $g_\mcC : \mcC \to \mcD$ such that $g_\mcB \circ f_\mcB = g_\mcC \circ f_\mcC$,

\item[(b)] {\em disjoint amalgamation property} if whenever $\mcA, \mcB \in \mbK$,
then there is $\mcC \in \mbK$ such that $C = A \cup B$, $\mcC \uhrc A = \mcA$ and $\mcC \uhrc B = \mcB$, and

\item[(c)] {\em free amalgamation property} if whenever $\mcA, \mcB \in \mbK$, then the free amalgam of $\mcA$
and $\mcB$ belongs to $\mbK$.
\end{itemize}
}\end{defin}

\begin{defin}\label{definition of age, permitted etc}{\rm
Let $\mcM$ be a $V$-structure, where the vocabulary $V$ contains only relation symbols.\\
(i) The {\em age} of $\mcM$ is the class of all finite $V$-structures that can be embedded into $\mcM$.\\
(ii) A finite $V$-structure is {\em permitted (with respect to $\mcM$)} if it belongs to the age of $\mcM$.
A finite $V$-structure is {\em forbidden (with respect to $\mcM$)} if it is not permitted.\\
(iii) A finite $V$-structure $\mcC$ is a {\em constraint (of $\mcM$)} if it is forbidden and every proper substructure of it is permitted.\\
(iv) $\mcM$ is {\em finitely constrained} (or {\em has only finitely many constraints}) if there are, up to isomorphism, 
only finitely many constraints of  $\mcM$.\\
(v) $\mcM$ is a {\em random structure} if it is homogeneous (implying that $V$ is finite) and for every $0 < k < \omega$,
every constraint of 
\[
\mcM \uhrc \{R \in V : \text{ the arity of $R$ is at most $k$}\}
\]
has cardinality at most $k$.
}\end{defin}

\noindent
The most well-known random structure is the {\em Rado graph} (often called {\em random graph} in model theory).
The Rado graph can be constructed in three interesting ways: as a Fra\"{i}ss\'{e} limit,  by a probabilistic construction on
finite graphs (via a so-called zero-one law), and by using a probability measure on the set of all graphs whose vertex set is the
natural numbers. (See for example \cite[Chapter~7.4]{Hod}.) Every random structure can be constructed in these three ways
(by straightforward generalizations of the corresponding procedures for graphs).

\subsection{Notions of triviality}\label{Notions of triviality}

If $\mcM$ is a structure then `$\acl_\mcM$' and `$\dcl_\mcM$' denote the algebraic closure and definable closure, respectively,
in $\mcM$.

\begin{defin}\label{triviality for acl}{\rm
Let $\mcM$ be a structure.\\
(i) We say that {\em $\acl_\mcM$ is trivial}, or that {\em algebraic closure (in $\mcM$) is trivial}, if
whenever $A \subseteq M$, $b \in M$ and $b \in \acl_\mcM(A)$, then $b \in \acl_\mcM(a)$ for some $a \in A$.\\
(ii) We say that {\em $\acl_\mcM$ is degenerate}, or that {\em algebraic closure (in $\mcM$) is degenerate}, if
for all $A \subseteq M$, $\acl_\mcM(A) = A$.
}\end{defin}

\begin{defin}\label{triviality of dependence}{\rm
Let $T$ be a simple theory.\\
(i) $T$ has {\em trivial dependence} if whenever $\mcM \models T$,
$A, B, C \subseteq M\meq$ and $A \underset{C}{\nind} B$, then $A \underset{C}{\nind} b$ for some $b \in B$. \\
(ii) Let $0 < k < \omega$.
$T$ has {\em $k$-trivial dependence for real elements over finite base sets} if whenever $\mcM  \models T$,
$A, B, C \subseteq M$ are finite and $A \underset{C}{\nind} B$, then $A \underset{C}{\nind} B'$ 
for some $B' \subseteq B$ with $|B'| \leq k$. \\
(iii) We say that a simple structure $\mcM$ has {\em trivial dependence}, or {\em
$k$-trivial dependence for real elements over finite base sets} if
its complete theory has the corresponding property.
}\end{defin}

\begin{observation}\label{observation about 2-transitivity and triviality}
Suppose that $\mcM$ is 2-transitive, simple and has 1-trivial dependence for real elements over finite base sets.
Then $\mcM$ is supersimple with SU-rank~1 and degenerate algebraic closure.
\end{observation}

\noindent
{\bf Proof.}
Suppose that $\mcN \models Th(\mcM)$, $B \subseteq N$ is finite, $a \in N \setminus B$  and $a \nind B$.
Since $\mcM$ has 1-trivial dependence for real elements over finite base sets, $a \nind b$ for some $b \in B$.
There is, by simplicity, $c \in N \setminus \{a\}$ (assuming that $\mcN$ is $\omega$-saturated) such that $a \ind c$.
Then we must have $\tp_\mcN(a, b) \neq \tp_\mcN(a, c)$ which contradicts 2-transitivity.
\hfill $\square$

\subsection{Imaginary elements, equivalence relations and $\omega$-categoricity}

We assume familiarity with the extension $\mcM\meq$ by imaginary elements of a structure $\mcM$.
Also recall that every infinite homogeneous structure is $\omega$-categorical, i.e. its complete theory is $\omega$-categorical.

\begin{fact}\label{basic facts about imaginaries}
Suppose that $\mcM$ is $\omega$-categorical.
\begin{itemize}
\item[(i)] If $B \subseteq M\meq$ is finite and $\bar{a} \in M\meq$, then $\tp_{\mcM\meq}(\bar{a} / \acl_{\mcM\meq}(B))$ is isolated.
\item[(ii)] If  $B \subseteq M\meq$ is finite and $p \in S_n^{\mcM\meq}(\acl_{\mcM\meq}(B))$ is realized in $\mcN\meq$ for some 
$\mcN \succcurlyeq \mcM$, then $p$ is realized in $\mcM\meq$; 
moreover, only finitely many types in $S_n^{\mcM\meq}(\acl_{\mcM\meq}(B))$ are realized by tuples from $M^n$.
\item[(iii)] Let $B \subseteq M$ be finite and let
$p_1, \ldots, p_k \in S_n^\mcM(B)$.
Then the following equivalence relation on $p_1(\mcM) \cup \ldots \cup p_k(\mcM)$ is $B$-definable in $\mcM$:
\[
\tp_{\mcM\meq}(\bar{x} / \acl_{\mcM\meq}(B)) = \tp_{\mcM\meq}(\bar{y} / \acl_{\mcM\meq}(B)). 
\]
\end{itemize}
\end{fact}

\noindent
Explanations of why~(i) and~(ii) above hold are given in \cite[Section 2.4]{AK}. 
Part~(iii) follows in a rather straightforward way from the previous parts.

\begin{fact}\label{facts about interpretability}
Suppose that $\mcN$ is interpretable in $\mcM$ using only finitely many parameters.
\begin{itemize}
\item[(i)] If $\mcM$ is $\omega$-categorical then so is $\mcN$.
\item[(ii)] If $\mcM$ is (super)simple then so is $\mcN$, and if the SU-rank of $\mcM$ is~1 
and $N \subseteq M$, then the SU-rank of $\mcN$ is~1.
\end{itemize}
\end{fact}

\noindent
Part~(i) is Theorem~7.3.8 in \cite{Hod}.
The claim about simplicity in~(ii) follows from Remarks~2.26 and~2.27 in \cite{Cas}.
The claim about supersimplicity follows from the fact that if some type of $\mcN$ divides over every finite subset of 
its set of parameters, then the same is true for some type of $\mcM\meq$ (so $\mcM$ would not be supersimple).
The final claim of~(ii) follows since (given the assumptions) 
every instance of dividing in $\mcN$ gives rise to a corresponding instance of dividing in $\mcM$.

\section{On definable structures}

\noindent
In this section we prove some technical results which will 
be used in Section~\ref{Extracting constraints from instances of dividing}
in the proof of Theorem~\ref{the rank is finite} via the use of 
Corollary~\ref{n-nontriviality implies contraint of size linear in n with sets on both sides and parameters}.

\begin{defin}\label{definition of structure M-p generated by p}{\rm
Let $\mcM$ be a $V$-structure where $V$ is a finite relational vocabulary and let $A \subseteq M$ be finite.
Furthermore, let $P = \{p_1, \ldots, p_n\}  \subseteq S_1^\mcM(A)$.
\begin{itemize}
\item[(i)] Let $V_A$ be a finite relational vocabulary such that $V \subseteq V_A$ and for every $R \in V$ of arity $r > 1$, every 
$0 < k < r$, every permutation $\pi$ of $\{1, \ldots, r\}$, 
and every $\bar{a} \in A^k$, $V_A$ has a relation symbol $Q_{R, \bar{a}, \pi}$ of arity $r-k$. 
We also assume that $V_A$ has no other symbols than those described.
Note that the maximal arity of $V_A$ is the same as the maximal arity of $V$.

\item[(ii)] Let $\mcM_P$ be the $V_A$-structure with universe $M_P = p_1(\mcM) \cup \ldots \cup p_n(\mcM)$ 
and where the symbols in $V_A$ are interpreted as follows:
\begin{itemize}
\item[(a)] If $R \in V$ has arity $r$, then $R^{\mcM_P} = R^\mcM \cap (M_P)^r$.

\item[(b)] If $Q_{R, \bar{a}, \pi} \in V_A \setminus V$ where $R \in V$ has arity $r$ and $|\bar{a}| = k$, then for every $\bar{b} \in (M_P)^{r-k}$,
$\bar{b} \in (Q_{R, \bar{a}, \pi})^{\mcM_P}$ if and only if $\pi(\bar{b}\bar{a}) \in R^\mcM$
(where $\pi(\bar{b}\bar{a}) = (c_{\pi(1)}, \ldots, c_{\pi(r)})$ if $\bar{b}\bar{a} = (c_1, \ldots, c_r)$).
\end{itemize}
\end{itemize}
}\end{defin}

\noindent
In the rest of this section we assume that $\mcM$ is a $V$-structure where $V$ is finite and relational, $A \subseteq M$ is finite and
$P = \{p_1, \ldots, p_n\}  \subseteq S_1^\mcM(A)$ is nonempty.
Furthermore, $V_A$ and $\mcM_P$ are as in Definition~\ref{definition of structure M-p generated by p}.
The following lemma is an immediate consequence of the definition of $\mcM_P$.

\begin{lem}\label{correspondence between quantifier free types}
Let $\bar{a}$ be an enumeration of $A$. For all $\bar{b}, \bar{b}' \in M_P$,
\[
\bar{b} \equiv_{\mcM_P}^\mrqf \bar{b}' \ \ \text{ if and only if } \ \ 
\bar{b}\bar{a} \equiv_\mcM^\mrqf \bar{b}'\bar{a}.
\]
\end{lem}

\noindent
For the rest of this section suppose that $\mcM$ is infinite and every type in $P$ is nonalgebraic.

\begin{lem}\label{M-p is simple and homogeneous}
If $\mcM$ is simple then $\mcM_P$ is simple. 
If $\mcM$ is homogeneous then $\mcM_P$ is homogeneous.
\end{lem}

\noindent
{\bf Proof.}
The first claim follows from Fact~\ref{facts about interpretability}.
Let $\bar{a}$ enumerate $A$.
Suppose that $\bar{b}, \bar{c} \in M_P$ and $\bar{b} \equiv_{\mcM_P}^\mrqf \bar{c}$.
Then Lemma~\ref{correspondence between quantifier free types} gives
$\bar{b}\bar{a} \equiv_\mcM^\mrqf \bar{c}\bar{a}$.
Since $\mcM$ is homogeneous there is an automorphism $f$ of $\mcM$ which takes $\bar{b}\bar{a}$ to $\bar{c}\bar{a}$.
Then $f$ sends every element of $M_P$ to an element of $M_P$.
By Lemma~\ref{correspondence between quantifier free types} again, the restriction of $f$ to $M_P$ is an
automorphism of $\mcM_P$. Thus $\mcM_P$ is homogeneous.
\hfill $\square$

\begin{lem}\label{independence in M and M-p}
Suppose that $\mcM$ is infinite and homogeneous.
For all $\bar{b}, \bar{c} \in M_P$, $\bar{b} \ind^{\mcM_P} \bar{c}$ if and only in $\bar{b} \underset{A}{\ind}^\mcM \bar{c}$.
\end{lem}

\noindent
{\bf Proof.}
Suppose that $\bar{b}, \bar{c} \in M_P$ and $\bar{b} \nind^{\mcM_P} \bar{c}$.
Then there are $\varphi(\bar{x}, \bar{y}) \in \tp_{\mcM_P}(\bar{b},\bar{c})$ and $\bar{c}_i$ in some elementary extension of $\mcM_P$, 
for $i < \omega$,
such that $\tp_{\mcM_P}(\bar{c}_i) = \tp_{\mcM_P}(\bar{c})$ for all $i$ and
$\{\varphi(\bar{x}, \bar{c}_i) : i < \omega\}$ is $k$-inconsistent for some $k < \omega$.
As $\mcM_P$ is homogeneous (by the previous lemma), we can assume
(by Fact~\ref{basic facts about imaginaries}) 
that all $\bar{c}_i$ belong to $\mcM_P$.
Without loss of generality we may assume that $\varphi(\bar{x}, \bar{y})$ isolates $\tp_{\mcM_P}(\bar{b},\bar{c})$ and
that $\varphi$ is quantifier-free.

Let $\bar{a}$ enumerate $A$.
By the homogeneity of $\mcM$ there is a quantifier-free $V$-formula $\psi(\bar{x}, \bar{y}, \bar{z})$ which isolates
$\tp_\mcM(\bar{b}, \bar{c}, \bar{a})$.
By Lemma~\ref{correspondence between quantifier free types} and since both $\mcM$ and $\mcM_P$ have elimination of quantifiers 
(by Lemma~\ref{M-p is simple and homogeneous}),
it follows that for all $\bar{b}', \bar{c}' \in M_P$,
\begin{equation}\label{going from psi to phi and conversely, first part}
\bar{b}'\bar{c}'\bar{a} \equiv_\mcM \bar{b}\bar{c}\bar{a} \ \Longleftrightarrow \ 
\bar{b}'\bar{c}' \equiv_{\mcM_P} \bar{b}\bar{c}
\end{equation}
and thus
\begin{equation}\label{going from psi to phi and conversely, second part}
\mcM \models \psi(\bar{b}', \bar{c}', \bar{a}) \ \Longleftrightarrow \ 
\mcM_P \models \varphi(\bar{b}', \bar{c}').
\end{equation}
From~(\ref{going from psi to phi and conversely, first part}) it follows that
\[
\bar{c}_i\bar{a} \ \equiv_\mcM \ \bar{c}\bar{a} \ \text{ for all $i < \omega$.}
\]
For a contradiction, assume that $\{\psi(\bar{x}, \bar{c}_i, \bar{a}) : i < \omega\}$ is not $k$-inconsistent with respect to $Th(\mcM)$.
Then (by Fact~\ref{basic facts about imaginaries}) there are $\bar{b}' \in M_P$ and $i_1, \ldots, i_k$ such that
\[
\mcM \models \bigwedge_{j=1}^k \psi(\bar{b}', \bar{c}_{i_j}, \bar{a}).
\]
By~(\ref{going from psi to phi and conversely, second part}) we get
\[
\mcM_P \models \bigwedge_{j=1}^k \varphi(\bar{b}', \bar{c}_{i_j})
\]
which implies that $\{\varphi(\bar{x}, \bar{c}_i) : i < \omega\}$ is not $k$-inconsistent, contradicting our assumption.
Thus we conclude that $\{\psi(\bar{x}, \bar{c}_i, \bar{a}) : i < \omega\}$ is $k$-inconsistent and it follows that
$\bar{b} \underset{A}{\nind}^\mcM \bar{c}$.

Now assume that $\bar{b} \underset{A}{\nind}^\mcM \bar{c}$ where $\bar{b}, \bar{c} \in M_P$.
Then there are a quantifier-free $\psi(\bar{x}, \bar{y}, \bar{z})$ such that 
$\psi(\bar{x}, \bar{c}, \bar{a})$ belongs to $\tp_\mcM(\bar{b}, \bar{c}, \bar{a})$ and isolates this type,
and $\bar{c}_i \in M_P$ such that $\tp_\mcM(\bar{c}_i, \bar{a}) = \tp_\mcM(\bar{c}, \bar{a})$ for all $i < \omega$
and $\{\psi(\bar{x}, \bar{c}_i, \bar{a}) : i < \omega\}$ is $k$-inconsistent for some $k < \omega$.
From~(\ref{going from psi to phi and conversely, first part}) it follows that
$\tp_{\mcM_P}(\bar{c}_i) = \tp_{\mcM_P}(\bar{c})$ for all $i$.
For the same reasons as when we proved the other direction, there is a quantifier-free $V_A$-formula $\varphi(\bar{x}, \bar{y})$
such that~(\ref{going from psi to phi and conversely, second part}) holds, 
and consequently $\varphi(\bar{x}, \bar{y})$ isolates $\tp_{\mcM_P}(\bar{b}, \bar{c})$.
For a contradiction, suppose that $\{\varphi(\bar{x}, \bar{c}_i) : i < \omega \}$ in not $k$-inconsistent.
Then there are $\bar{b}' \in M_P$ and $i_1, \ldots, i_k$ such that 
$\mcM_P \models \bigwedge_{j=1}^k \varphi(\bar{b}', \bar{c}_{i_j})$. 
By~(\ref{going from psi to phi and conversely, second part}) we get
$\mcM \models \bigwedge_{j=1}^k \psi(\bar{b}', \bar{c}_{i_j}, \bar{a})$, 
which implies that $\{\psi(\bar{b}', \bar{c}_i, \bar{a}) : i < \omega\}$ is not $k$-inconsistent, contradicting our assumption.
Hence we conclude that $\bar{b} \nind^{\mcM_P} \bar{c}$.
\hfill $\square$

\begin{lem}\label{for every constraint of M-p there is a larger constraint of $M}
Suppose that $\mcM$ is homogeneous.
For every constraint $\mcC$ of $\mcM_P$ there is a constraint $\mcC^*$ of $\mcM$ such that $|C| \leq |C^*|$.
\end{lem}

\noindent
{\bf Proof.}
Let $\mcC$ be a constraint of $\mcM_P$.
Let $\mcD$ be the $V$-structure (where $V$ is the vocabulary of $\mcM$) with universe $D = C \cup A$ and satisfying the following conditions:
\begin{itemize}
\item[(a)] For every $R \in V$ and every $\bar{c} \in C$ of appropriate length, $\bar{c} \in R^\mcD$ if and only if $\bar{c} \in R^\mcC$.

\item[(b)] For every $R \in V$ and every $\bar{a} \in A$ of appropriate length,
$\bar{a} \in R^\mcD$ if and only if $\bar{a} \in R^\mcM$.

\item[(c)] For every $R \in V$ of arity $r > 1$, every $0 < k < r$, every permutation $\pi$ of $\{1, \ldots, r\}$,
every $\bar{c} \in C^{r-k}$ and every $\bar{a} \in A^k$,
$\pi(\bar{c}\bar{a}) \in R^\mcD$ if and only if $\bar{c} \in (Q_{R,\bar{a}, \pi})^{\mcC}$
(where $Q_{R,\bar{a}, \pi}$ is like in Definition~\ref{definition of structure M-p generated by p}).
\end{itemize}

\medskip

\noindent
{\bf Claim.} $\mcD$ is forbidden with respect to $\mcM$.

\medskip
\noindent
{\em Proof of the claim.}
Suppose that $\mcD$ is not forbidden with respect to $\mcM$.
Then, by the homogeneity of $\mcM$, there is $\mcD' \subseteq \mcM$ and an isomorphism $f$ from $\mcD$ to $\mcD'$.
Let $C' = f(C)$ and $A' = f(A)$. Then the restriction of $f$ to $A$ is an isomorphism from $\mcD \uhrc A$ to $\mcM \uhrc A'$.
So by homogeneity again, there is $C'' \subseteq M$ and an isomorphism from $\mcD$ to $\mcM \uhrc A \cup C''$ which extends $f \uhrc A$.
But then, by the definition of $\mcM_P$, $\mcM_P \uhrc C'' \cong \mcC$ which contradicts that $\mcC$ is a constraint of $\mcM_P$.
\hfill $\square$

\medskip
\noindent
{\bf Claim.} There is $B \subseteq A$ such that either
\begin{itemize}
\item[(i)] $B = \es$ and $\mcD \uhrc C$ is forbidden with respect to $\mcM$, or
\item[(ii)] $B \neq \es$, $\mcD \uhrc B \cup C$ is forbidden with respect to $\mcM$
and for every $b \in B$, $\mcD \uhrc (B \cup C) \setminus \{b\}$ is permitted with respect to $\mcM$.
\end{itemize}

\medskip
\noindent
{\em Proof of the claim.}
If for all $a \in A$, $\mcD \uhrc D \setminus \{a\}$ is permitted with respect to $\mcM$,
then take $B = A$ and we are done.
Otherwise there is some $a_1 \in A$ such that $\mcD \uhrc D \setminus \{a_1\}$ is forbidden with respect to $\mcM$.
If for all $a \in A \setminus \{a_1\}$, $\mcD \uhrc D \setminus \{a_1, a\}$ is permitted with respect to $\mcM$,
then take $B = A \setminus \{a_1\}$ and we are done.
Otherwise there is some $a_2 \in A \setminus \{a_1\}$ such that 
$\mcD \uhrc D \setminus \{a_1, a_2\}$ is forbidden with respect to $\mcM$.
By continuing in this way we eventually find $B \subseteq A$ such that~(i) or~(ii) holds.
\hfill $\square$

\medskip
\noindent
Let $B \subseteq A$ satisfy~(i) or~(ii) of the last claim.
In order to show that $\mcD \uhrc B \cup C$ is a constraint of $\mcM$ it suffices to show that if $c \in C$,
then $\mcD \uhrc (B \cup C) \setminus \{c\}$ is permitted with respect to $\mcM$.
So let $c \in C$.
As $\mcC$ is a constraint of $\mcM_P$, $\mcC \uhrc C \setminus \{c\}$ is permitted with respect to $\mcM_P$.
As $\mcM_P$ is homogeneous there is an embedding $f$ of $\mcC \uhrc C \setminus \{c\}$ into $\mcM_P$.
From the definitions of $\mcM_P$ and $\mcD$ it follows that if $f$ is extended to $(B \cup C) \setminus \{c\}$ 
in such a way that the extension fixes
all elements of $B$ pointwise, then this extension is an embedding of $\mcD \uhrc (B \cup C) \setminus \{c\}$ into $\mcM$,
so $\mcD \uhrc (B \cup C) \setminus \{c\}$ is permitted.
Now we have proved that $\mcC^* = \mcD \uhrc (B \cup C)$ is a constraint of $\mcM$ and clearly
$|C| \leq |C^*|$ since $C \subseteq C^*$.
\hfill $\square$

\section{Extracting constraints from instances of dividing}\label{Extracting constraints from instances of dividing}

\noindent
In this section we prove our first main results, Theorems~\ref{the rank is finite} and~\ref{trivial acl if the rank is one}.
As indicated in the title of this section, the core of the proofs is to extract a constraint from every
instance of dividing. The main technical result in this section is 
Proposition~\ref{large parameter set in a dividing formula implies large constraint},
from which we get a number of corollaries which extend its usefulness.

\begin{theor}\label{the rank is finite}
Suppose that $\mcM$ is ternary, homogeneous and simple.
If $\mcM$ has only finitely many constraints, then it is supersimple with finite SU-rank;
moreover, if $n$ is the cardinality of the largest constraint and $k = 1 + n \cdot |S_3^\mcM(\es)|$, then 
$\mcM$ has $k$-trivial dependence for real elements over finite base sets.
\end{theor}

\begin{theor}\label{trivial acl if the rank is one}
Suppose that $\mcM$ is ternary, homogeneous and supersimple with SU-rank~1.
If $\mcM$ has only finitely many constraints, then $\acl_\mcM$ is trivial.
\end{theor}

\noindent
Note that if $\mcM$ in Theorem~\ref{trivial acl if the rank is one} is, in addition, 2-transitive then $\acl_\mcM$ is degenerate.
Also observe that the conclusion of Theorem~\ref{trivial acl if the rank is one}
implies that $\mcM$ has 1-trivial dependence for real elements over finite base sets;
by (for example) the argument in Remark~6.6 in \cite{KopPrimitive} it follows that $\mcM$
has trivial dependence.
The proofs of Theorems~\ref{the rank is finite} and~\ref{trivial acl if the rank is one}
are given after Lemma~\ref{uniformity of k in k-inconsistency} below.

\begin{rem}\label{remark on main theorems}{\rm
Theorems~\ref{the rank is finite} and~\ref{trivial acl if the rank is one} can be strengthened in the following ways,
as proved in~\cite[Corollaries~3.4--3.5]{Kop17b} which I wrote after this article but it was refereed and 
accepted for publication before this article.
First, the conclusion of Theorem~\ref{the rank is finite} can be strengthened to say that $\mcM$ has trivial dependence.
Secondly, in Theorem~\ref{trivial acl if the rank is one} the assumption that $\mcM$ has only finitely many constraints can be removed.
}\end{rem}

\begin{prop}\label{large parameter set in a dividing formula implies large constraint}
Let $\mcM$ be an infinite ternary homogeneous structure.
Suppose that $a, b \in M$, $E, F \subseteq M$ are finite, some formula
$\varphi(x, b, \bar{e}) \in \tp_\mcM(a / bE)$
divides over $EF$, and 
for every proper subset $E' \subset E$, every formula $\psi(x, b, \bar{e}') \in \tp_\mcM(a / bE')$ does not divide over $EF$.
Then $\mcM$ has a constraint with at least $3 + |E|/|S_3^\mcM(\es)|$ elements.
\end{prop}

\noindent
{\bf Proof.}
Let $a, b \in M$, $E, F \subseteq M$ and 
$\varphi(x, b\bar{e}) \in \tp_\mcM(a / bE)$
be such that the assumptions of the proposition are satisfied.
We rename these elements and sets as follows:
$a_n = a$, $a_{n-1} = b$, $A = E = \{a_1, \ldots, a_{n-2}\}$.
Furthermore, let $p(x_1, \ldots, x_n) = \tp_\mcM(a_1, \ldots, a_n)$ and when convenient we identity
$p$ notationally with a quantifier-free formula that isolates it.
By assumption, 
$p(a_1, \ldots, a_{n-1}, x_n)$ divides over $AF$.
So there is an $AF$-indiscernible sequence $(b_i : i < \omega)$ such that
$b_0 = a_{n-1}$ and for some $1 < k < \omega$, 
\[
\{ p(a_1, \ldots, a_{n-2}, b_i, x_n) : \ i < \omega \} \ \text{ is $k$-inconsistent, but not $l$-inconsistent if $l < k$.}
\]
For every $j \in \{1, \ldots, n-2\}$, let
\begin{align*}
& A_j = A \setminus \{a_j\}, \text{ and} \\
&\text{let $p_j$ be the restriction of $p$ to the variables in $\{x_1, \ldots, x_n\} \setminus \{x_j\}$.}
\end{align*}
When convenient we notationally identify $p_j$ with a quantifier-free formula that isolates it.
By assumption, for every $j = 1, \ldots, n-2$, 
\[
p_j(a_1, \ldots, a_{j-1}, a_{j+1}, \ldots, a_{n-2}, a_{n-1}, x_n) \ \text{ does not divide over $AF$,}
\]
and therefore, for every $j = 1, \ldots, n-2$,
\[
\{ p_j(a_1, \ldots, a_{j-1}, a_{j+1}, \ldots, a_{n-2}, b_i, x_n) : i < \omega \} \ \text{ is $m$-consistent for every $m < \omega$.}
\]
By Ramsey's theorem, for every $j = 1, \ldots, n-2$ there are $a_n^j \in M$ and distinct $b_1^j, \ldots, b_k^j \in \{b_i : i < \omega\}$ such that
\begin{align}\label{property of b-l-j}
&\mcM \models \bigwedge_{l = 1}^k p_j(a_1, \ldots, a_{j-1}, a_{j+1}, \ldots, a_{n-2}, b_l^j, a_n^j) \ \text{ and}\\
&(a_n^j, b_l^j, b_{l'}^j) \equiv_\mcM (a_n^j, b_s^j, b_{s'}^j) \ \text{ for all $1 \leq l < l' \leq k$ and all $1 \leq s < s' \leq k$.} \nonumber
\end{align}
Since $\mcM$ is ternary with elimination of quantifiers it follows that 
\begin{equation}\label{b-1-j is a-n-j indiscernible}
\text{$(b_1^j, \ldots, b_k^j)$ is an $\{a_n^j\}$-indiscernible sequence for each $j$.}
\end{equation}
By elimination of quantifiers there are only finitely many 3-types over $\es$. 
Hence there is a number $\tau$, at most as large as the number of 3-types (of distinct elements) over $\es$,
and a partition $A = X_1 \cup \ldots \cup X_\tau$ such that $X_1, \ldots, X_\tau$
are the equivalence classes of the following equivalence relation on $A$:
\[
a_j \sim a_{j'} \ \Longleftrightarrow \ (a_n^j, b_1^j, b_2^j) \equiv_\mcM (a_n^{j'}, b_1^{j'}, b_2^{j'}).
\]
Let $a, b'_1, \ldots, b'_k$ be distinct elements. (It does not matter where they come from.)
For each $m \in \{1, \ldots, \tau\}$ let $\mcB_m$ be the structure with universe
\[
B_m \ = \ A \cup \{b'_1, \ldots, b'_k, a\}
\]
such that 
\begin{align*}
&\tp_{\mcB_m}^\mrqf(a_1, \ldots, a_{n-2}, b'_1, \ldots, b'_k) \ = \ 
\tp_\mcM^\mrqf(a_1, \ldots, a_{n-2}, b_1, \ldots, b_k), \\
&\tp_{\mcB_m}^\mrqf(a_1, \ldots, a_{n-2}, b'_s, a)  \ = \ 
\tp_\mcM^\mrqf(a_1, \ldots, a_{n-2}, b_1, a_n)  \ \text{ for all } s = 1, \ldots, k, \text { and} \\
&\tp_{\mcB_m}^\mrqf(b'_1, \ldots, b'_k, a) \ = \ 
\tp_\mcM^\mrqf(b_1^j, \ldots, b_k^j, a_n^j), \ \text{ where $j$ is any number in } \\
&\text{$\{1, \ldots, n-2\}$ such that $a_j  \in X_m$. (Recall that $\tp_\mcM^\mrqf(b_1^j, \ldots, b_k^j, a_n^j)$}\\
&\text{only depends on the $\sim$-class of $a_j$.)}
\end{align*}
Observe that if $m \neq m'$, then the interpretations of relation symbols in $\mcB_m$ and in $\mcB_{m'}$
only differ on triples containing $a$ and two elements from $\{b'_1, \ldots, b'_k\}$.
Also note that each $\mcB_m$ is forbidden, due to the
$k$-inconsistency of $\{ p(a_1, \ldots, a_{n-2}, b_i, x_n) : \ i < \omega \}$.

Some $X_m$ must have cardinality at least $(n-2)/\tau$.
Without loss of generality, assume that the cardinality of $X_1$ is at least $(n-2)/\tau$.
From~(\ref{property of b-l-j}) it follows that
\begin{equation}\label{removing an element from X-0}
\text{for every $c \in X_1$, \ $\mcB_1 \uhrc (A \setminus \{c\}) \cup \{b'_1, \ldots, b'_k, a\}$ is permitted.}
\end{equation}
Since $A \cup \{b_i : i < \omega\} \subseteq M$ it follows from the construction of $\mcB_1$ that
\begin{equation}\label{removing a}
\mcB_1 \uhrc A \cup \{b'_1, \ldots, b'_k\} \ \text{ is permitted}.
\end{equation}
By the choice of the elements $b_i$, $i < \omega$, and the construction of $\mcB_1$ it follows that
\begin{equation}\label{removing all b'-l but one}
\text{for each $l = 1, \ldots, k$, \ $\mcB_1 \uhrc A \cup \{b'_l\} \cup \{a\}$ is permitted.}
\end{equation}

\medskip
\noindent
{\bf Claim.} There is $C \subseteq A \setminus X_1$ such that 
\begin{align}\label{removing c from C}
&\mcB_1 \uhrc X_1 \cup C \cup \{b'_1, \ldots, b'_k, a\} \text{ is forbidden and}\\ 
&\text{for every $c \in C$, }
\mcB_1 \uhrc X_1 \cup (C \setminus \{c\}) \cup \{b'_1, \ldots, b'_k, a\} \text{ is permitted.} \nonumber
\end{align}

\medskip
\noindent
{\em Proof of the claim.}
If $\mcB_1 \uhrc X_1 \cup \{b'_1, \ldots, b'_k, a\}$ is forbidden, then,
by taking $C = \es$,~(\ref{removing c from C}) holds for trivial reasons.
Suppose that $\mcB_1 \uhrc X_1 \cup \{b'_1, \ldots, b'_k, a\}$ is permitted.
If there is $c_1 \in A \setminus X_1$ such that 
$\mcB_1 \uhrc X_1 \cup \{c_1, b'_1, \ldots, b'_k, a\}$ is forbidden then take $C = \{c_1\}$ and we are done.
Otherwise $\mcB_1 \uhrc X_1 \cup \{c, b'_1, \ldots, b'_k, a\}$ is permitted for all $c \in A \setminus X_1$.
If there are $c_1, c_2 \in A \setminus X_1$ such that 
$\mcB_1 \uhrc X_1 \cup \{c_1, c_2, b'_1, \ldots, b'_k, a\}$ is forbidden then take $C = \{c_1, c_2\}$ and we are done.
Otherwise $\mcB_1 \uhrc X_1 \cup \{c_1, c_2, b'_1, \ldots, b'_k, a\}$ is permitted for all $c_1, c_2 \in A \setminus X_1$.
Since $\mcB_1$ is forbidden we will eventually, if we continue in this way, find $C  \subseteq A \setminus X_1$ such 
that~(\ref{removing c from C}) holds.
\hfill $\square$

\medskip

\noindent
According to the claim there is $C  \subseteq A \setminus X_1$ such that~(\ref{removing c from C}) holds.
This together with~(\ref{removing an element from X-0}) implies that
\begin{equation}\label{removing c from the union of X-0 and C}
\text{for every $c \in X_1 \cup C$, \ $\mcB_1 \uhrc \big((X_1 \cup C) \setminus \{c\}\big) \cup \{b'_1, \ldots, b'_k, a\}$ is permitted.}
\end{equation}

\medskip

\noindent
{\bf Claim.} There is $D \subseteq \{b'_1, \ldots, b'_k\}$ such that $|D| \geq 2$,
\begin{align}\label{removing an element from D}
&\mcB_1 \uhrc X_1 \cup C \cup D \cup \{a\} \ \text{ is forbidden and}\\
&\text{for every } d \in D, \ \mcB_1 \uhrc X_1 \cup C \cup (D \setminus \{d\}) \cup \{a\} \ \text{ is permitted.} \nonumber
\end{align}

\medskip

\noindent
{\em Proof of the claim.}
If for all $d \in \{b'_1, \ldots, b'_k\}$, $\mcB_1 \uhrc X_1 \cup C \cup \big(\{b'_1, \ldots, b'_k\} \setminus \{d\}\big) \cup \{a\}$
is permitted, then take $D = \{b'_1, \ldots, b'_k\}$ and we are done.
Otherwise there is $d_1 \in \{b'_1, \ldots, b'_k\}$ such that 
$\mcB_1 \uhrc X_1 \cup C \cup \big(\{b'_1, \ldots, b'_k\} \setminus \{d_1\}\big) \cup \{a\}$ is forbidden.
If for all $d \in \{b'_1, \ldots, b'_k\} \setminus \{d_1\}$,
$\mcB_1 \uhrc X_1 \cup C \cup \big(\{b'_1, \ldots, b'_k\} \setminus \{d_1, d\}\big) \cup \{a\}$ is permitted, 
then take $D = \{b'_1, \ldots, b'_k\} \setminus \{d_1\}$ and we are done.
Otherwise there is $d_2 \in \{b'_1, \ldots, b'_k\} \setminus \{d_1\}$ such that 
$\mcB_1 \uhrc X_1 \cup C \cup \big(\{b'_1, \ldots, b'_k\} \setminus \{d_1, d_2\}\big) \cup \{a\}$ is forbidden.
From~(\ref{removing all b'-l but one}) it follows that by continuing in this way we will eventually find
$D \subseteq \{b'_1, \ldots, b'_k\}$ with at least two elements such that~(\ref{removing an element from D}) holds.
\hfill $\square$

\medskip
\noindent
By the claim let $D \subseteq \{b'_1, \ldots, b'_k\}$ have at least two elements and satisfy~(\ref{removing an element from D}).
It now follows from~(\ref{removing a}),~(\ref{removing c from the union of X-0 and C}) 
and~(\ref{removing an element from D}) that $\mcB_1 \uhrc X_1 \cup C \cup D \cup \{a\}$ is a constraint.
Moreover, by the choice of $X_1$ and $D$, $|X_1| \geq (n-2)/\tau$
and $|D| \geq 2$, where $\tau$ is at most as large as the number of 3-types over $\es$.
\hfill $\square$

\begin{cor}\label{n-nontriviality implies contraint of size linear in n}
Let $\mcM$ be an infinite ternary homogeneous structure.
Suppose that $n \geq 3$, $a_1, \ldots, a_n \in M$,
$\tp_\mcM(a_n / a_{n-1}, \ldots, a_1)$ divides over $\{a_{n-2}, \ldots, a_1\}$, but
for every proper subset $A' \subset \{a_{n-2}, \ldots, a_1\}$,
$\tp_\mcM(a_n / a_{n-1}A')$ does not divide over $A'$.
Then $\mcM$ has a constraint with at least $3 + (n-2)/|S_3^\mcM(\es)|$ elements.
\end{cor}

\noindent
{\bf Proof.}
Suppose that $n \geq 3$ and that $a_1, \ldots, a_n \in M$ satisfy the assumptions of the corollary.
Let $A = \{a_1, \ldots, a_{n-2}\}$, $F = \es$ and $p(x_1, \ldots, x_n) = \tp_\mcM(a_1, \ldots, a_n)$.
Now we can argue exactly as in the proof of Proposition~\ref{large parameter set in a dividing formula implies large constraint}.
(However, I do not see how to use the statement of Proposition~\ref{large parameter set in a dividing formula implies large constraint}
directly to get Corollary~\ref{n-nontriviality implies contraint of size linear in n}.)
\hfill $\square$

\begin{cor}\label{corollary to n-nontriviality implies contraint of size linear in n}
Let $\mcM$ be a ternary simple homogeneous structure.
Suppose that $a \in M$, $A, B \subset M$ are finite sets, $|B| \geq 2$, $a \underset{A}{\nind} B$ 
and $a \underset{A}{\ind} B'$ for every proper subset $B' \subset B$.
Then $\mcM$ has a constraint with at least $3 + (|B| - 1)/|S_3^\mcM(\es)|$ elements.
\end{cor}

\noindent
{\bf Proof.}
Let $\mcM$, $a \in M$ and $A, B \subset M$ satisfy the assumptions of the corollary.
Note that the assumptions imply that $\acl_\mcM(A) \cap B = \es$.
Let $p_1, \ldots, p_m$ enumerate all types in $S_1^\mcM(A)$ which are realized in $\{a\} \cup B$
and let $P = \{p_1, \ldots, p_m\}$, so all types in $P$ are nonalgebraic.
Furthermore, let $\mcM_P$ be as in Definition~\ref{definition of structure M-p generated by p}.
By Lemmas~\ref{M-p is simple and homogeneous}
and~\ref{independence in M and M-p},
$\mcM_P$ is homogeneous, simple, $a \nind^{\mcM_P} B$ and $a \ind^{\mcM_P} B'$ for every proper subset $B' \subset B$.
Let the elements of $B$ be enumerated as $a_1, \ldots, a_{n-1}$ (so $n \geq 3$ as $|B| \geq 2$) and let $a_n = a$.
Then $\tp_{\mcM_P}(a_n / a_{n-1}, \ldots, a_1)$ divides over $\es$ and, for every proper subset
$B' \subset \{a_{n-1}, \ldots, a_1\}$, $\tp_{\mcM_P}(a_n / B')$ does not divide over $\es$.
By transitivity, $\tp_{\mcM_P}(a_n / a_{n-1}, \ldots, a_1)$ divides over $\{a_{n-2}, \ldots, a_1\}$.
By monotonicity, for every proper $A' \subset \{a_{n-2}, \ldots, a_1\}$, $\tp_{\mcM_P}(a_n / a_{n-1}A')$
does not divide over $A'$.
By Corollary~\ref{n-nontriviality implies contraint of size linear in n},
$\mcM_P$ has a constraint with at least $3 + (|B|- 1)/|S_3^\mcM(\es)|$ elements.
By Lemma~\ref{for every constraint of M-p there is a larger constraint of $M},
so does $\mcM$.
\hfill $\square$

\begin{cor}\label{n-nontriviality implies contraint of size linear in n with sets on both sides}
Let $\mcM$ be a ternary simple homogeneous structure.
Suppose that $A, B \subseteq M$ are finite, $A \nind B$, $A' \ind B$ for every proper subset $A' \subset A$,
and $A \ind B'$ for every proper subset $B' \subset B$.
If $|B| \geq 2$ then
$\mcM$ has a constraint with at least $3 + (|B|-1)/|S_3^\mcM(\es)|$ elements.
\end{cor}

\noindent
{\bf Proof.}
Let $A$ and $B$ satisfy the assumptions of the corollary,
so in particular $A \cap B = \es$.
Let $A = \{a_1, \ldots, a_r\}$.
Then $\{a_2, \ldots, a_r\} \ind B$ and by symmetry, $B \ind \{a_2, \ldots, a_r\}$ (and $B \nind A$),
so we get $B \underset{\{a_2, \ldots, a_r\}}{\nind} a_1$ by transitivity, and then
$a_1 \underset{\{a_2, \ldots, a_r\}}{\nind} B$ by symmetry.
If $a_1 \underset{\{a_2, \ldots, a_r\}}{\nind} B'$ for some proper subset $B' \subset B$,
then (by symmetry and monotonicity) $\{a_1, \ldots, a_r\} \nind B'$, which contradicts our assumption.
Hence $a_1 \underset {\{a_2, \ldots, a_r\}}{\ind} B'$ for every proper subset $B' \subset B$.
Now the result follows from 
Corollary~\ref{corollary to n-nontriviality implies contraint of size linear in n}.
\hfill $\square$

\begin{cor}\label{n-nontriviality implies contraint of size linear in n with sets on both sides and parameters}
Let $\mcM$ be a ternary simple homogeneous structure.
Suppose that $A, B, C$ $\subseteq M$ are finite, $A \underset{C}{\nind} B$ and $A' \underset{C}{\ind} B$ for every proper subset $A' \subset A$,
and $A \underset{C}{\ind} B'$ for every proper subset $B' \subset B$.
Then $\mcM$ has a constraint with at least $3 + (|B| - 1)/|S_3^\mcM(\es)|$ elements.
\end{cor}

\noindent
{\bf Proof.}
Like the proof of Corollary~\ref{corollary to n-nontriviality implies contraint of size linear in n},
but we use Corollary~\ref{n-nontriviality implies contraint of size linear in n with sets on both sides}
instead of 
Corollary~\ref{n-nontriviality implies contraint of size linear in n}.
\hfill $\square$
\\

\noindent
By \cite[Proposition~2.5.4]{Kim_book}, the conclusion of the next lemma is equivalent to saying that $Th(\mcM)$ is 
a {\em low} theory.

\begin{lem}\label{uniformity of k in k-inconsistency}
Let $\mcM$ be a ternary simple homogeneous structure.
For every formula $\varphi(x, \bar{y})$ without parameters, there is $k_\varphi < \omega$,
depending only on $\varphi$, such that if $\bar{c}_i$, $i < \omega$, is an indiscernible sequence
and $\{\varphi(x, \bar{c}_i) :  i < \omega\}$ is inconsistent, then it is $k_\varphi$-inconsistent.
\end{lem}

\noindent
{\bf Proof.}
Suppose that $(\bar{c}_i : i < \omega)$ is an indiscernible sequence.
Since $\mcM$ is ternary with elimination of quantifiers, for every $m < \omega$ and
all $i_1 < \ldots < i_m < \omega$, $\tp_\mcM(\bar{c}_{i_1}, \ldots, \bar{c}_{i_m})$ is determined by
$\tp_\mcM(\bar{c}_0, \bar{c}_1, \bar{c}_2)$.
Let us call $\tp_\mcM(\bar{c}_0, \bar{c}_1, \bar{c}_2)$ the type of the indiscernible sequence $(\bar{c}_i : i < \omega)$.
By $\omega$-categoricity, for every $s < \omega$, there are only finitely many types of indiscernible sequences
$(\bar{c}_i : i < \omega)$, such that $|\bar{c}_0| = s$.
By homogeneity, if $(\bar{c}_i : i < \omega)$ and $(\bar{d}_i : i < \omega)$ are indiscernible sequences of the same type
and $\{\varphi(x, \bar{c}_i) : i < \omega\}$ is $k$-consistent, then so is 
$\{\varphi(x, \bar{d}_i) : i < \omega\}$.
\hfill $\square$
\\

\noindent
\subsection{Proof of Theorem~\ref{the rank is finite}}

Let $\mcM$ be ternary, homogeneous and simple with only finitely many constraints.
Let $\rho$ be the cardinality of a constraint of $\mcM$ of maximal cardinality.
Note that by Proposition~\ref{large parameter set in a dividing formula implies large constraint},
whenever the assumptions of this proposition are satisfied, then $3 + |E|\big/|S_3^\mcM(\es)| \leq \rho$,
so $|E|  \leq (\rho - 3) \cdot |S_3^\mcM(\es)|$.

To prove that $\mcM$ is supersimple with finite SU-rank it essentially suffices to prove the following lemma.

\begin{lem}\label{arbitrarily long dividing chains do not exist}
There is $n_0 < \omega$ such that there do not exist $n_0 < n < \omega$, $a \in M$ and 
finite subsets $B_0 \subset B_1 \subset \ldots \subset B_n \subset M$ such that $a \underset{B_i}{\nind} B_{i+1}$
for all $i < n$.
\end{lem}

\noindent
{\bf Proof.}
Suppose that $a \in M$, $B_0 \subset B_1 \subset \ldots B_n$ are finite subsets of $M$
and $a \underset{B_i}{\nind} B_{i+1}$ for every $i < n$.
For every $i < n$, we can choose a formula $\varphi_{i+1}(x, \bar{b}_{i+1}) \in \tp_\mcM(a, B_{i+1})$ such that
$\varphi_{i+1}(x, \bar{b}_{i+1})$ divides over $B_i$ and for every
$\bar{b}$ such that $\rng(\bar{b})$ is a proper subset of $\rng(\bar{b}_{i+1})$,
every formula $\psi(x, \bar{b}) \in \tp_\mcM(a, \bar{b})$ does not divide over $B_i$.
By Proposition~\ref{large parameter set in a dividing formula implies large constraint}
we have $|\bar{b}_i| < \rho \cdot |S_3^\mcM(\es)|$ for all $0 < i \leq n$.
By $\omega$-categoricity, there are, up to equivalence in $\mcM$, only finitely many, say $t$,
formulas with at most $\rho \cdot |S_3^\mcM(\es)| + 1$ free variables. 
Note that the numbers $\rho$, $|S_3^\mcM(\es)|$ and $t$ depend only on $Th(\mcM)$.
Moreover, by Lemma~\ref{uniformity of k in k-inconsistency}, 
there is $k < \omega$ such that, for every $i < n$, 
$\varphi_{i+1}(x, \bar{b}_{i+1})$ $k$-divides over $B_i$.\footnote{
That $\varphi_{i+1}(x, \bar{b}_{i+1})$ $k$-divides over $B_i$ means that there is a $B_i$-indiscernible sequence
$(\bar{c}_\alpha : \alpha < \omega)$ such that $\bar{b}_{i+1} = \bar{c}_0$ and
$\{\varphi_{i+1}(x, \bar{c}_\alpha) : \alpha < \omega\}$ is $k$-inconsistent.}
Again, $k$ depends only on $Th(\mcM)$.
Let $f(x) = \lfloor x/t \rfloor$, so $f$ depends only on $Th(\mcM)$ and not on $n$.
Then, if $n$ is sufficiently large, there are $0 < m_0 < \ldots < m_{f(n)} \leq n$ such that
$\varphi_{m_i}(x, \bar{y})$ is equivalent, in $\mcM$, to $\varphi_{m_j}(x, \bar{y})$ for all $i < j \leq f(n)$ and 
$\varphi_{m_i}(x, \bar{b}_{m_i})$ $k$-divides over $B_{m_i - 1}$ for all $i \leq f(n)$.
Let us rename $\varphi_{m_0}(x, \bar{y})$ by $\psi_n(x, \bar{y})$.
From the definition of $k$-dividing it follows that, for every $i < f(n)$,
$\psi_n(x, \bar{b}_{m_i})$ $k$-divides over $\bigcup_{j < i}\rng(\bar{b}_{m_j})$.
Consequently we have
$D(\tp_\mcM(a / B_0), \psi_n(x, \bar{y}), k) \geq f(n)$, so $D(x=x, \psi_n(x, \bar{y}), k) \geq f(n)$ 
(by \cite[Lemma~2.3.4]{Kim_book} for example).

As $\rho, |S_3^\mcM(\es)|, t, k$ and the function $f(x)$ only depend on $Th(\mcM)$
and $\psi_n$ is always a formula in at most $\rho \cdot |S_3^\mcM(\es)| + 1$ free variables
(of which there are only $t$ choices up to equivalence modulo $Th(\mcM)$),
it follows that if such $a \in M$ and finite subsets $B_0 \subset B_1 \subset \ldots \subset B_n$ of $M$ exist
for every $n < \omega$, then for some $\psi(x, \bar{y})$ with at most 
$\rho \cdot |S_3^\mcM(\es)| + 1$ free variables, $D(x=x, \psi(x, \bar{y}), k) \geq n$ for every $n < \omega$, so 
$D(x=x, \psi(x, \bar{y}), k)$ is infinite.
This would imply, using \cite[Proposition~2.3.7]{Kim_book}, that $\mcM$ is not simple.
Hence, there is $m < \omega$ such that if $n > m$ then there do not exist
$a \in M$ and finite subsets $B_0 \subset B_1 \subset \ldots B_n$ of $M$
such that $a \underset{B_i}{\nind} B_{i+1}$ for every $i < n$.
\hfill $\square$
\\

\noindent
By the finite character of dividing and 
lemma~\ref{arbitrarily long dividing chains do not exist}
it follows that $\mcM$ is supersimple with finite SU-rank.
The second part of the theorem is a direct consequence of 
Corollary~\ref{n-nontriviality implies contraint of size linear in n with sets on both sides and parameters}.
This concludes the proof of Theorem~\ref{the rank is finite}

\subsection{Proof of Theorem~\ref{trivial acl if the rank is one}}

Suppose that $\mcM$ is ternary, homogeneous, supersimple with SU-rank~1 and with only finitely many constraints.
Then, as the SU-rank of $\mcM$ is~1, for all $a \in M$ and $B, C \subseteq M$, $a \underset{C}{\nind} B$ if and only if
$a \in \acl_\mcM(B \cup C) \setminus \acl_\mcM(C)$.
Suppose for a contradiction, that $\acl_\mcM$ is not trivial. 
Then there is finite $A = \{a_1, \ldots, a_n\} \subset M$ such that $n \geq 3$, $A$ is not independent but every proper subset of $A$ independent.
Since $a_1 \ind \{a_2, \ldots, a_{n-1}\}$ we find, using the existence of nondividing extensions,
$a'_2, \ldots, a'_{n-1}$ such that $\tp_\mcM(a_1, a'_2, \ldots, a'_{n-1}) = \tp_\mcM(a_1, a_2, \ldots, a_{n-1})$
and 
\[
\{a'_2, \ldots, a'_{n-1}\} \underset{a_1}{\ind} \{a_2, \ldots, a_{n-1}\}.
\]
Then there is $a'_n \in \acl_\mcM(a_1, a'_2, \ldots, a'_{n-2})$ such that 
\[
\tp_\mcM(a_1, a'_2, \ldots, a'_n) = \tp_\mcM(a_1, a_2, \ldots, a_n).
\]
Now it is straightforward to verify that $\{a_2, \ldots, a_n, a'_2, \ldots, a'_n\}$ is not independent but every
proper subset of it is independent.
By repeating the procedure we obtain arbitrarily large finite dependent $B \subset M$ such that every proper subset of $B$ 
is independent.
Now Corollary~\ref{n-nontriviality implies contraint of size linear in n} implies that 
$\mcM$ has arbitrarily large constraints, contradicting the assumption that $\mcM$ is finitely constrained.\footnote{
One can also use Lemma~3 in \cite{Goode} the proof of which works out in the present context, but I chose to give a simpler proof here.
}

\section{Definable equivalence relations and weakly isolated contraints}\label{Definable equivalence relations and weakly isolated contraints}

\noindent
In order to understand the fine structure of a simple structure $\mcM$ it is useful to understand
the fine structure on definable subsets of $M\meq$ (and in particular of $M$) of SU-rank~1.
Also, by Observation~\ref{observation about 2-transitivity and triviality}, if $\mcM$ is 2-transitive
and has trivial dependence, then its SU-rank is~1 and $\acl_\mcM$ is trivial.
Moreover, there are 
(see Section~\ref{Uncountably many ternary homogeneous simple structures}), 
for a vocabulary $V$ with only one ternary relation symbol, uncountably many homogeneous
2-transitive supersimple $V$-structures with SU-rank~1 and degenerate algebraic closure.
Thus there is reason for dealing with the special case when the SU-rank is~1, even under the assumption that the algebraic closure is degenerate.
The ultimate goal would be to find some sort of classification of such structures.
We do not reach that far here.
But we do reveal some general connections between
the nature of the constraints of a structure
and the existence of nontrivial definable (with parameters)
equivalence relations in it.
Here, an equivalence relation is called {\em nontrivial} if it has more than one equivalence class and has at least
one equivalence class with more than one element.
The question whether there is a ternary, homogeneous, supersimple structure $\mcM$ with SU-rank~1
that contains a finite subset $A \subset M$  such that $|A| \geq 3$, $A$ is not independent (over $\es$) but every proper
subset of $A$ is independent remains open. 
Note that by Theorem~\ref{trivial acl if the rank is one} such an example, if it exists, must have infinitely many constraints.

\begin{defin}\label{definition of weakly isolated constraint}{\rm
Let $V$ be a finite relational vocabulary.\\
(i) Let $\mcA$ and $\mcB$ be $V$-structures with the same universe and let $a, b, c$ be distinct elements from their
universe. We say that $\mcA$ and $\mcB$ are {\em $(a, b, c)$-neighbours} if for every finite sequence $\bar{d}$ of
elements from their universe such that $\{a, b, c\} \not\subseteq \rng(\bar{d})$, 
$\tp_\mcA^\mrqf(\bar{d}) = \tp_\mcB^\mrqf(\bar{d})$.\\
(ii) Suppose that $\mcM$ is a ternary homogeneous structure and let $\mcC$ be a constraint of $\mcM$.
We say that $\mcC$ is {\em weakly isolated
(with respect to $\mcM$, or the age of $\mcM$)} 
if for every choice of distinct $a, b, c \in C$, $\mcC$ has a $(a, b, c)$-neighbour
which is permitted with respect to $\mcM$.
We say that $\mcC$ is {\em isolated
(with respect to $\mcM$, or the age of $\mcM$)} 
if for every choice of distinct $a, b, c \in C$, every $(a, b, c)$-neighbour of $\mcC$ is
permitted with respect to $\mcM$.
}\end{defin}

\noindent
All examples known to the author of ternary 2-transitive homogeneous supersimple structures with SU-rank~1 have only
weakly isolated constraints.
Note that if we work with 3-hypergraphs, then the notions `isolated' and `weakly isolated' coincide.
It will turn out that the existence of constraints which are {\em not} weakly isolated is related to the existence
of nontrivial definable (with parameters) equivalence relations.
But first we consider the case when the age of a structure has the free amalgamation property.

\begin{prop}\label{free amalg implies no nontrivial equivalence relation}
Suppose that $\mcM$ is homogeneous and that its age has the free amalgamation property.
Then for every finite $A \subset M$ and every nonalgebraic $p(x) \in S_1^\mcM(A)$,
there is no nontrivial $A$-definable equivalence relation on $p(\mcM)$.
\end{prop}

\noindent
{\bf Proof.}
Let $\mcM$ satisfy the assumptions of the proposition. 
(We can assume that $\mcM$ is infinite, which would follow if we assumed that $\mcM$ was simple.
For if $\mcM$ is finite then all types are algebraic and the result follows automatically.)
It implies that $\acl_\mcM$ is degenerate,
for if $\bar{b} \in M$ and $a \in M \setminus \rng(\bar{b})$, then, by free amalgamation
(in fact disjoint amalgamation suffices), for every $n$ there are distinct $a_1, \ldots, a_n \in M$ such that
$a_i\bar{b} \equiv_\mcM^\mrqf a\bar{b}$ for all $i$ and by elimination of quantifiers we get
$a_i\bar{b} \equiv_\mcM a\bar{b}$ for all $i$, so $\tp_\mcM(a / \bar{b})$ cannot be algebraic.
Suppose for a contradiction that $A \subset M$ is finite, $p \in S_1^\mcM(A)$ is nonalgebraic and that 
$E$ is a nontrivial $A$-definable equivalence relation on $p(\mcM)$.
Since $\acl_\mcM$ is degenerate it follows that every equivalence class of $E$ is infinite.
Since $\mcM$ is homogeneous it follows that $E$ is defined by a quantifier-free formula with parameters from $A$.

First suppose that if $a, b \in p(\mcM)$ are distinct and $E(a, b)$, then there is some $\bar{c} \in A$ such that some
permutation of $ab\bar{c}$ is a relationship.
As $E$ is nontrivial on $p(\mcM)$ there are distinct $a, b \in p(\mcM)$ such that $E(a, b)$.
Let $\mcA_0 = \mcM \uhrc (A \cup \{a\})$ and $\mcA_1 = \mcM \uhrc (A \cup \{a, b\})$.
Furthermore, let $\mcA_2$ be an isomorphic copy of $\mcA_1$ such that $\mcA_0 \subseteq \mcA_2$ and 
$A_2 \setminus (A \cup \{a\}) = \{b'\}$ where $b' \neq b$.
By free amalgamation the free amalgam, say $\mcB$, of $\mcA_1$ and $\mcA_2$ over $\mcA_0$ is permitted.
Without loss of generality we may assume that $\mcB \subseteq \mcM$, so that also $b'$ belongs to $M$.
As $E$ is an equivalence relation and $E(a, b)$ and $E(a, b')$ we must have $E(b, b')$
But as $\mcB$ is the free amalgam of $\mcA_1$ and $\mcA_2$,
there is no $\bar{c} \in A$ such that some permutation of $bb'\bar{c}$ is a relationship.
Since $E(b, b')$ holds this contradicts our assumption.

Now suppose that there are distinct $a, b \in p(\mcM)$ such that $E(a, b)$ and there is no $\bar{c} \in A$ such that
some permutation of $ab\bar{c}$ is a relationship.
Since $\mcM$ has elimination of quantifiers (and the vocabulary is finite and relational) it follows that
\begin{align}\label{E holds if abc does not form a relationship}
&\text{if $a, b \in p(\mcM)$ are distinct and there is no $\bar{c} \in A$ such that some}\\
&\text{permutation of $ab\bar{c}$ is a relationship, then $E(a, b)$.} \nonumber
\end{align}
Let $a_1, a_2 \in p(\mcM)$ be distinct but otherwise arbitrary.
We will show that $E(a_1, a_2)$ which contradicts that $E$ is nontrivial.
Let $\mcA_0 = \mcM \uhrc (A \cup \{a_1\})$ and $\mcA_1 = \mcM \uhrc (A \cup \{a_1, a_2\})$.
Since every $E$-class is infinite there is $b \in p(\mcM) \setminus \{a_1\}$ such that $E(a_1, b)$.
Let $\mcA_2 = \mcM \uhrc (A \cup \{a_1, b\})$. 
By free amalgamation the free amalgam, say $\mcB$, of $\mcA_1$ and $\mcA_2$ over $\mcA_0$ is permitted.
Without loss of generality we may assume that $\mcB \subseteq \mcM$.
It follows from~(\ref{E holds if abc does not form a relationship}) that $E(a_2, b)$.
Since $E$ is an equivalence relation and $E(a_1, b)$ we get $E(a_1, a_2)$.
\hfill $\square$
\\

\noindent
The next two lemmas are just straightforward observations, most likely noticed by others before,
but we nevertheless give their proofs.

\begin{lem}\label{nondegenerate acl and equivalence relations}
Suppose that $\mcM$ is $\omega$-categorical and supersimple with SU-rank~1.
Also assume that there are finite $A \subseteq M$ and distinct $b, c \in M \setminus \acl_\mcM(A)$ such that
$c \in \acl_\mcM(bA) \setminus \acl_\mcM(A)$. 
Then, letting $p_1 = \tp_\mcM(b / A)$ and $p_2 = \tp_\mcM(c / A)$, there is a 
nontrivial $A$-definable equivalence relation on $p_1(\mcM) \cup p_2(\mcM)$
with only finite equivalence classes.
\end{lem}

\noindent
{\bf Proof.}
Let $A \subseteq M$ be finite and let $b, c \in M \setminus \acl_\mcM(A)$ be distinct such that
$c \in \acl_\mcM(bA) \setminus \acl_\mcM(A)$. 
Let $p_1(x) = \tp_\mcM(b / A)$ and $p_2(x) = \tp_\mcM(c / A)$.
As $\mcM$ has SU-rank~1 
(so $(M, \acl_\mcM)$ is a pregeometry) and is $\omega$-categorical
it follows that if we define
$E(x,y)$ if and only if $x \in \acl_\mcM(yA)$, then $E$ is a nontrivial 
$A$-definable equivalence relation on $p_1(\mcM) \cup p_2(\mcM)$
with only finite equivalence classes.
\hfill $\square$

\begin{lem}\label{degenerate acl, types over imaginaries and equivalence relations}
Suppose that $\mcM$ is $\omega$-categorical.
The following are equivalent for any finite $A \subset M$:
\begin{itemize}
\item[(a)] There are $b, c \in M \setminus \acl_\mcM(A)$ such that 
$\tp_\mcM(b / A) = \tp_\mcM(c / A)$ and \\
$\tp_{\mcM\meq}(b / \acl_{\mcM\meq}(A)) \neq \tp_{\mcM\meq}(c / \acl_{\mcM\meq}(A))$.
\item[(b)] There are a nonalgebraic $p \in S_1^\mcM(A)$ and a nontrivial
$A$-definable equivalence relation on $p(\mcM)$ with only finitely many equivalence classes all of which are infinite.
\end{itemize}
\end{lem}

\noindent
{\bf Proof.}
The implication from~(a) to~(b) follows from the fact (Fact~\ref{basic facts about imaginaries}~(iii)) that the equivalence relation
\[
\tp_{\mcM\meq}(x / \acl_{\mcM\meq}(A)) = \tp_{\mcM\meq}(y / \acl_{\mcM\meq}(A))
\]
on $p(\mcM)$ is $A$-definable and has only finitely many equivalence classes 
all of which are infinite (since $p = tp_\mcM(b / A)$ is nonalgebraic).

Suppose that~(b) holds and let $E$ be a nontrivial $A$-definable equivalence relation on $p(\mcM)$.
Let $\bar{a}$ enumerate $A$ and let $n = |\bar{a}|$.
Define the following equivalence relation on $M^{n+1}$:
$E'(x\bar{x}, y\bar{y})$ if and only if 
\[
\bar{x} = \bar{y} \ \text{ and } \ \tp_{\mcM\meq}(x / \acl_{\mcM\meq}(\bar{x})) = \tp_{\mcM\meq}(y / \acl_{\mcM\meq}(\bar{x})).
\]
By a straightforward argument, using, Fact~\ref{basic facts about imaginaries}, it follows that $E'$ is a $\es$-definable equivalence relation.
As we assume that $E$ is nontrivial on $p(\mcM)$ there are $b, c \in p(\mcM)$ such that $\neg E(b, c)$, from which it
follows (since $E$ is $A$-definable) that $\neg E'(b\bar{a}, c\bar{a})$.
By Fact~\ref{basic facts about imaginaries}, only finitely many complete types over 
$\acl_{\mcM\meq}(A)$ are realized in $M$, so the $E'$-classes of $b\bar{a}$ and of $c\bar{a}$ belong
to $\acl_{\mcM\meq}(A)$. Consequently $b$ and $c$ have different types over
$\acl_{\mcM\meq}(A)$.
\hfill $\square$
\\

\noindent
We have seen above that if the algebraic closure is degenerate, 
then the existence of definable equivalence relations is related to the
existence of elements having (for some finite set $A$) the same type over $A$ but different types over $\acl_{\mcM\meq}(A)$.
The next result and its corollary implies that if $\mcM$ is finitely constrained,
then every nontrivial definable equivalence relation 
(on a nonalgebraic type) 
is determined by equivalence relations which are definable over sets of parameters
of bounded size.

\begin{theor}\label{definable equivalence relations with many parameters}
Let $\mcM$ be ternary, homogeneous and supersimple with SU-rank 1 and degenerate algebraic closure.
Suppose that $A \subset M$ is finite, that $p \in S_1^\mcM(A)$ is nonalgebraic and that there 
are $b, c  \in p(\mcM)$ such that 
$\tp_{\mcM\meq}(b / \acl_{\mcM\meq}(A)) \neq 
\tp_{\mcM\meq}(c / \acl_{\mcM\meq}(A))$
and for every proper subset $A' \subset A$,
$\tp_{\mcM\meq}(b / \acl_{\mcM\meq}(A')) = \tp_{\mcM\meq}(c / \acl_{\mcM\meq}(A'))$.
Then $\mcM$ has a constraint $\mcC$ with at least $3 + |A|/|S_3^\mcM(\es)|$ elements.
\end{theor}

\noindent
{\bf Proof.}
Suppose that $A \subseteq M$ is finite, 
let $p \in S_1^\mcM(A)$ be a nonalgebraic type
and suppose that there are $b, c \in p(\mcM)$ are such that
$\tp_{\mcM\meq}(b / \acl_{\mcM\meq}(A)) \neq 
\tp_{\mcM\meq}(c / \acl_{\mcM\meq}(A))$ and for every 
proper subset $A' \subset A$,
$\tp_{\mcM\meq}(b / \acl_{\mcM\meq}(A')) = \tp_{\mcM\meq}(c / \acl_{\mcM\meq}(A'))$.
If we define $E$ on $p(\mcM)$ by
$E(x, y)$ if and only if 
$\tp_{\mcM\meq}(x / \acl_{\mcM\meq}(A)) =
\tp_{\mcM\meq}(y / \acl_{\mcM\meq}(A))$,
then (by Fact~\ref{basic facts about imaginaries})
$E$ is an $A$-definable equivalence relation with only finitely many equivalence classes all of which are infinite
(since $\acl_\mcM$ is degenerate).
By the choice of $b$ and $c$, $\neg E(b, c)$.
As all equivalence classes are infinite there are
$a, a' \in p(\mcM) \setminus \{b, c\}$ such that $E(a, b)$ and $E(a', c)$.

Let $A = \{a_1, \ldots, a_n\}$, $\bar{a} = (a_1, \ldots, a_n)$ and let
$\mcF$ be a structure with universe $F = A \cup \{a, b, c\}$ and such that
\begin{align}\label{definition of F}
&\tp_\mcF^\mrqf(b, c, \bar{a}) = \tp_\mcM^\mrqf(b, c, \bar{a}), \\
&\tp_\mcF^\mrqf(a, b, \bar{a}) = \tp_\mcM^\mrqf(a, b, \bar{a}), \ \text{ and} \nonumber\\
&\tp_\mcF^\mrqf(a, c, \bar{a}) = \tp_\mcM^\mrqf(a', c, \bar{a}). \nonumber
\end{align}
Since $\mcM$ has elimination of quantifiers there is a quantifier-free formula $\varphi(x, y, \bar{z})$ 
without parameters such that $\varphi(x, y, \bar{a})$ defines $E$.
We claim that $\mcF$ is forbidden.
If not, then there is an embedding $f : \mcF \to \mcM$.
By quantifier-elimination, $\tp_\mcM(f(\bar{a})) = \tp_\mcM(\bar{a})$ and consequenctly
$\varphi(x, y, f(\bar{a}))$ defines an equivalence relation on $q(\mcM)$ where $q(x) = \{\varphi(x, f(\bar{a}')) : \varphi(x, \bar{a}') \in p\}$.
But as $f$ is an embedding and $\varphi$ is quantifier-free we get
\[
\mcM \models \varphi(f(a), f(b), f(\bar{a})) \wedge \varphi(f(a), f(c), f(\bar{a})) \wedge \neg\varphi(b, c, f(\bar{a}))
\]
which contradicts that $\varphi(x, y, f(\bar{a}))$ defines an equivalence relation on $q(\mcM)$.
Thus $\mcF$ is forbidden, no matter how $\tp_\mcF^\mrqf(a, b, c)$ is chosen.
Hence every $(a, b, c)$-neighbour of $\mcF$ is forbidden.

\medskip

\noindent
{\bf Claim.} For every $i \in \{1, \ldots, n\}$, $\mcF \uhrc (F \setminus \{a_i\})$ has an $(a, b, c)$-neighbour $\mcP_i$ which is permitted.

\medskip

\noindent
{\bf Proof.}
Let $i \in \{1, \ldots, n\}$ and let $A_i = A \setminus \{a_i\}$.
By the choice of $b$ and $c$ and since $E(a, b)$ and $E(a', c)$ we have
$\tp_{\mcM\meq}(a /\acl_{\mcM\meq}(A_i)) = \tp_{\mcM\meq}(a' /\acl_{\mcM\meq}(A_i))$.
Since we assume that $\mcM$ is supersimple with SU-rank 1 and $\acl_\mcM$ is degenerate we have
$a \underset{A_i}{\ind} b$, $a' \underset{A_i}{\ind} c$ and $b \underset{A_i}{\ind} c$.
The independence theorem of simple theories (and the fact that $\mcM$ is $\omega$-saturated)
thus implies that there is $a'' \in M$ such that
\[
\tp_\mcM(a'' / bA_i) = \tp_\mcM(a / bA_i) \ \text{ and }  \ \tp_\mcM(a'' / cA_i) = \tp_\mcM(a' / cA_i).
\]
Let $\mcP_i$ be  structure with universe $A_i \cup \{a, b, c\}$
such that the map which sends $a''$ to $a$ and every element in $A_i \cup \{b, c\}$ to itself
is an isomorphism from $\mcM \uhrc (A_i \cup \{a'', b, c\})$ to $\mcP_i$. 
By the choice of $a''$ and since the language is ternary it follows that $\mcP_i$ is permitted and is an $(a, b, c)$-neighbour of $\mcF$.
\hfill $\square$

\medskip

\noindent
Let $\tau = |S_3^\mcM(\es)|$.
By the pigeon hole principle there are $m \geq n/\tau$ and $B \subseteq A$ such that $|B| = m$ and for all
$a_i, a_j \in B$, $\tp_{\mcP_i}^\mrqf(a, b, c) = \tp_{\mcP_j}^\mrqf(a, b, c)$.
To simplify notation we assume that $B = \{a_1, \ldots, a_m\}$.
Now we have that 
\[
\mcP_i \uhrc (A \setminus \{a_i, a_j\}) \cup \{a, b, c\} \ = \ \mcP_j \uhrc (A \setminus \{a_i, a_j\}) \cup \{a, b, c\} \ 
\text{ whenever } 1 \leq i \leq j \leq m.
\]
Therefore we can define a structure $\mcF'$ with universe $F'  = A \cup \{a, b, c\}$ such that
\begin{align*}
&\tp_{\mcF'}^\mrqf(b, c, a_1, \ldots, a_n) = \tp_\mcM^\mrqf(b, c, a_1, \ldots, a_n) \  \text{ and, for all } i = 1, \ldots, m, \\
&\tp_{\mcF'}^\mrqf(a, b, c, a_1, \ldots, a_{i-1}, a_{i+1}, \ldots, a_n) \ = \
\tp_{\mcP_i}^\mrqf(a, b, c, a_1, \ldots, a_{i-1}, a_{i+1}, \ldots, a_n).
\end{align*}
Note that for all triples $(d_1, d_2, d_3) \in (F')^3$ except permutations of $(a, b, c)$, we have
$\tp_{\mcF'}^\mrqf(d_1, d_2, d_3) = \tp_\mcF^\mrqf(d_1, d_2, d_3)$.
Therefore $\mcF'$ is forbidden.
By construction we also have that
\begin{equation}\label{F' minus an element is permitted}
\text{for every } d \in \{a, b, c, a_1, \ldots, a_m\}, \ \mcF' \uhrc F' \setminus \{d\} \ \text{ is permitted.}
\end{equation}
The next claim is proved similarly as the first claim in the proof of 
Proposition~\ref{large parameter set in a dividing formula implies large constraint}
so we omit its proof.

\medskip
\noindent
{\bf Claim.} There is $D \subseteq A \setminus B$ such that $\mcF' \uhrc B \cup D \cup \{a, b, c\}$ is forbidden and for every
$d \in D$, $\mcP \uhrc (B \cup D \cup \{a, b, c\}) \setminus \{d\}$ is permitted.

\medskip
\noindent
Let $D \subseteq A \setminus B$ be as in the claim.
From~(\ref{F' minus an element is permitted}) it follows that $\mcF' \uhrc B \cup D \cup \{a, b, c\}$ is a constraint
with at least $3 + |B| \geq 3 + |A|/|S_3^\mcM(\es)|$ elements.
\hfill $\square$

\begin{cor}\label{equivalence relations and finitely many constraints}
Suppose that $\mcM$ is ternary, homogeneous and supersimple with SU-rank 1 and degenerate algebraic closure.
Furthermore, assume that $\mcM$ has only finitely many constraints.
Then there is $n < \omega$ such that for every finite $A \subseteq M$, every nonalgebraic $p \in S_1^\mcM(A)$
and all $a, a' \in p(\mcM)$, if 
$\tp_{\mcM\meq}(a / \acl_{\mcM\meq}(B)) = \tp_{\mcM\meq}(a' / \acl_{\mcM\meq}(B))$ for all 
$B \subseteq A$ with $|B| \leq n$, then
$\tp_{\mcM\meq}(a / \acl_{\mcM\meq}(A)) = \tp_{\mcM\meq}(a' / \acl_{\mcM\meq}(A))$.
\end{cor}

\noindent
{\bf Proof.}
Let $\mcM$ satisfy the assumptions of the corollary.
If the conclusion is false, then, for every $n < \omega$, there is finite $A_n \subset M$ such that $|A_n| > n$,
and nonalgebraic $p \in S_1^\mcM(A_n)$ and $b, c \in p(\mcM)$ such that 
$\tp_{\mcM\meq}(b / \acl_{\mcM\meq}(A_n)) \neq \tp_{\mcM\meq}(c / \acl_{\mcM\meq}(A_n))$
and for every proper subset $B \subset A_n$ such that $|B| \leq n$,
$\tp_{\mcM\meq}(b / \acl_{\mcM\meq}(B)) = \tp_{\mcM\meq}(c / \acl_{\mcM\meq}(B))$.
By considering subsets of $A_n$ of cardinality at least $n$ we find, for each $n$, $A'_n \subseteq A_n$ such that 
$n < |A'_n|$,
$\tp_{\mcM\meq}(b / \acl_{\mcM\meq}(A'_n)) \neq \tp_{\mcM\meq}(c / \acl_{\mcM\meq}(A'_n))$ and
for every proper subset $B \subset A'_n$, 
$\tp_{\mcM\meq}(b / \acl_{\mcM\meq}(B)) = \tp_{\mcM\meq}(c / \acl_{\mcM\meq}(B))$.
By Theorem~\ref{definable equivalence relations with many parameters},
there is a constraint with at least $3 + n/|S_3^\mcM(\es)|$ elements for each $n < \omega$.
This contradicts that $\mcM$ has only finitely many constraints.
\hfill $\square$
\\

\noindent
The next theorem and its corollaries relates weakly isolated constraints to the existence of definable equivalence relations.

\begin{theor}\label{weak isolation and nontrivial equivalence relations}
Suppose that $\mcM$ is ternary, homogeneous and supersimple with SU-rank 1.
If $\mcC$ is a constraint of $\mcM$ with at least three elements, then at least one of the following 
two conditions holds:
\begin{itemize}
\item[(a)] $\mcC$ is weakly isolated.
\item[(b)] There is $\mcD \subseteq \mcC$ with $|D| = |C| - 3$ such that for every embedding $f : \mcD \to \mcM$,
\begin{itemize}
\item[(i)] there are (not necessarily distinct) nonalgebraic types $p_1, p_2 \in S_1^\mcM(f(D))$ and a 
nontrivial $f(D)$-definable equivalence relation on $p_1(\mcM) \cup p_2(\mcM)$ with only finite equivalence classes, or
\item[(ii)] there is a nonalgebraic type $p \in S_1^\mcM(f(D))$ and a nontrivial $f(D)$-definable equivalence relation on $p(\mcM)$
with only finitely many equivalence classes all of which are infinite.
\end{itemize}
\end{itemize}
\end{theor}

\noindent
{\bf Proof.}
Let $\mcC$ be a constraint of $\mcM$ with at least three elements, let $a, b, c \in C$ be distinct
and let $D = C \setminus \{a, b, c\}$.
Since $\mcC \uhrc D \cup \{b, c\}$ is permitted and is $\mcM$ homogeneous we may, without loss of generality, assume that
$\mcC \uhrc D \cup \{b, c\} \subseteq \mcM$.
Since $\mcC \uhrc D \{a, b\}$ is permitted we may also assume that 
$\mcC \uhrc D \cup \{a, b\} \subseteq \mcM$.
As $\mcC \uhrc D \cup \{a, c\}$ is permitted there is $a' \in M$ such that if $\bar{d}$ enumerates $D$,
then $\tp_\mcM^\mrqf(a', c, \bar{d}) = \tp_\mcC^\mrqf(a, c, \bar{d})$.
Note that by elimination of quantifiers in $\mcM$ we have $\tp_\mcM(a' / D) = \tp_\mcM(a / D)$.

First suppose that
\begin{align}
&b \underset{D}{\ind} c, \label{b is independent from c over D} \\
&a \underset{D}{\ind} b, \label{a is independent from b over D} \\
&a' \underset{D}{\ind} c, \ \text{ and} \label{a' is independent from c over D} \\
&\tp_{\mcM\meq}(a / \acl_{\mcM\meq}(D)) = \tp_{\mcM\meq}(a' / \acl_{\mcM\meq}(D)).
\label{a and a' have the same type in Meq over D}
\end{align}
Then the independence theorem for simple theories implies that there is $a'' \in M$
such that
\begin{equation}\label{properties of a''}
\tp_\mcM(a'' / bD) = \tp_\mcM(a / bD), \ \tp_\mcM(a'' / cD) = \tp_\mcM(a' / cD), \ \text{ and } \ a'' \underset{D}{\ind} bc.
\end{equation}
Let $\mcC^*$ be the unique structure with universe $C$ such that the following map $f$ is an isomorphism from 
$\mcM \uhrc D \cup \{a'', b, c\}$ to $\mcC^*$: $f(a'') = a$ and $f$ is the identity on $D \cup \{b, c\}$.
Then $\mcC^*$ is permitted and, since the language is ternary, an $(a, b, c)$-neighbour of $\mcC$.
If~(\ref{b is independent from c over D})--(\ref{a and a' have the same type in Meq over D}) 
hold for every choice of distinct $a, b, c \in C$, then $\mcC$ is weakly isolated.

We conclude that if $\mcC$ is not weakly isolated, then there are distinct $a, b , c \in C$ such that at least one of
(\ref{b is independent from c over D})--(\ref{a and a' have the same type in Meq over D}) fails.
If~(\ref{b is independent from c over D}) fails, then, since the SU-rank of $\mcM$ is 1, $b \in \acl_\mcM(cD) \setminus \acl_\mcM(D)$.
It follows from 
Lemma~\ref{nondegenerate acl and equivalence relations} 
that if $p_1 = \tp_\mcM(b)$ and $p_2 = \tp_\mcM(c)$, 
then there is a nontrivial $D$-definable equivalence relation on $p_1(\mcM) \cup p_2(\mcM)$ with only finite classes.
If~(\ref{a is independent from b over D}) or~(\ref{a' is independent from c over D})
fails then we argue in the same way.
If~(\ref{b is independent from c over D})--(\ref{a' is independent from c over D}) hold and
~(\ref{a and a' have the same type in Meq over D}) fails, then, 
by Lemma~\ref{degenerate acl, types over imaginaries and equivalence relations}, 
with $p = \tp_\mcM(a)$ there is a nontrivial $D$-definable equivalence relation on $p(\mcM)$ with only finitely
many equivalence classes all of which are infinite.
\hfill $\square$

\begin{cor}\label{no nontrivial equivalence relations imply only weakly isolated constraints}
Suppose that $\mcM$ is ternary, homogeneous and supersimple with SU-rank~1.
Furthermore, suppose that for every finite $A \subset M$ and all (not necessarily distinct) nonalgebraic $p_1, p_2 \in S_1^\mcM(A)$
there is no nontrivial $A$-definable equivalence relation on $p_1(\mcM) \cup p_2(\mcM)$ except, if $p_1 \neq p_2$, for the 
equivalence relation: ``$x$ and $y$ have the same complete type over $A$''.
Then every constraint of $\mcM$ is weakly isolated.
\end{cor}

\noindent
{\bf Proof.}
If the assumptions of the corollary are satisfied, then, for every constraint $\mcC$, condition~(b) of 
Theorem~\ref{weak isolation and nontrivial equivalence relations}
fails, so $\mcC$ must be weakly isolated by the same theorem.
\hfill $\square$

\begin{cor}\label{free amalg implies on weakly isolated constraints}
Suppose that $\mcM$ is ternary, homogeneous, supersimple with SU-rank~1 and that its age has the free amalgamation property.
Then every constraint of $\mcM$ is weakly isolated.
\end{cor}

\noindent
{\bf Proof.}
Note that if $\mcM$ satisfies the assumptions then $\acl_\mcM$ is degenerate.
Therefore condition~(i) of Theorem~\ref{weak isolation and nontrivial equivalence relations} cannot hold.
The conclusion now follows from 
Proposition~\ref{free amalg implies no nontrivial equivalence relation}
and Theorem~\ref{weak isolation and nontrivial equivalence relations}.
\hfill $\square$
\\

\noindent
For some questions left unanswered, see Section~\ref{Problems}.

\section{Definability of binary random structures}\label{Definability of binary random structures}

\noindent
Suppose that $\mcM$ is ternary, homogeneous, supersimple with SU-rank~1 and degenerate algebraic closure.
We saw in the previous section that if for every finite $A \subset M$ and every nonalgebraic type $p \in S_1^\mcM(A)$,
there is no nontrivial $A$-definable equivalence relation on $p(\mcM)$,
then every constraint of $\mcM$ is weakly isolated. 
The following question remains open: Is there $\mcM$ with the properties assumed above such that 
there is no nontrivial $\es$-definable equivalence relation on $M$, but for some finite $A \subset M$
and some nonalgebraic $p \in S_1^\mcM(A)$, there is a nontrivial $A$-definable equivalence relation on $p(\mcM)$?
Proposition~\ref{strong primitivity over a singleton} below implies
that if $\mcM$ is, in addition,
symmetric and 2-transitive, and $A$ is a {\em singleton}, 
then there is no nontrivial $A$-definable equivalence relation on $M \setminus A$.
A direct consequence (Corollary~\ref{M-p-minus is a binary random structure}) is that, with the same assumptions on $\mcM$,
a binary random structure is definable in $\mcM$ using only one parameter.
Thus, if $\mcM$ is {\em not} 3-transitive, then the Rado graph is definable in $\mcM$.
The first result of this section has a similar conclusion, but starts from more general assumptions.
It also gives information about the fine structure of $\mcM$.

\begin{defin}\label{definition of extension property}{\rm
Let $\mcM$ be a {\em binary} structure and $E$ an equivalence relation on $M$.
Here we say that $\mcM$ {\em satisfies extension properties relative to $E$}
if whenever $0 < n < \omega$, $a_1, \ldots, a_n, b_1, \ldots, b_n \in M$, all $b_1, \ldots, b_n$ belong to the same
$E$-class and $b_i \neq a_i$ for all $i$, then there is $b \in M$ such that
$\tp_\mcM(a_i, b) = \tp_\mcM(a_i, b_i)$ for all $i = 1, \ldots, n$.
}\end{defin}

\noindent
Although the context here is partly the same as that studied by Ahlman \cite{Ahl} 
(in particular Theorem~5.7 of \cite{Ahl}) the definition
of extension properties is different here, because we can do with the simpler definition above.
(The more complicated definition of `$\xi$-extension properties' in \cite{Ahl} makes sense in a more general context, including that
of considering a class of finite structures and asymptotic probabilities.)
Extension properties relative to some equivalence relation can be useful for understanding what the age of the structure looks like.
Moreover, under suitable circumstances extension properties can be used to carry out a back-and-forth argument which establishes
that a structure has certain properties, like being homogeneous.

\begin{defin}\label{definition of M-p-minus}{\rm
Let $\mcM$ be an $\omega$-categorical $V$-structure where $V$ is a finite relational vocabulary.
Suppose that $A \subseteq M$ is finite, $P \subseteq  S_1^\mcM(A)$ and let $V_A$ be the vocabulary and $\mcM_P$ the $V_A$-structure in
Definition~\ref{definition of structure M-p generated by p}.
Then $\mcM_P^-$  be the reduct of $\mcM_P$ to the vocabulary $V_A \setminus V$.
}\end{defin}

\noindent
Observe that if $\mcM$ is ternary then $\mcM_P^-$ is binary.

\begin{prop}\label{M-p-minus is omega-categorical}
Suppose that $\mcM$ is ternary, homogeneous, supersimple with SU-rank 1 and with degenerate $\acl_\mcM$.
Let $A \subset M$ be finite and let $P \subseteq S_1^\mcM(A)$ be a nonempty set of nonalgebraic types.
\begin{itemize}
\item[(i)] $\mcM_P^-$ is supersimple with SU-rank 1.
\item[(ii)] The algebraic closure operator in $\mcM_P^-$ is degenerate. 
\item[(iii)] $\mcM_P^-$ satisfies extension properties relative to the following $\es$-definable (in $\mcM_P^-$) equivalence relation:
$E^P(x,y)$ if and only if 
\[
\tp_{(\mcM_P^-)\meq}(x / \acl_{(\mcM_P^-)\meq}(\es)) = 
\tp_{(\mcM_P^-)\meq}(y / \acl_{(\mcM_P^-)\meq}(\es)).
\]
\item[(iv)] $E^P$ is definable (in $\mcM_P^-$) by a quantifier-free formula without parameters.
\end{itemize}
\end{prop}

\noindent
{\bf Proof.}
We first note that $\mcM_P^-$ is interpretable in $\mcM$ with finitely many parameters (those in $A$).
So by Fact~\ref{facts about interpretability}, $\mcM_P^-$ is $\omega$-categorical and supersimple with SU-rank~1, so~(i) is proved.
Part (ii) is a straightforward consequence of the fact that $\mcM_P^-$ is interpretable in $\mcM$
in such a way that every element $a \in \mcM_P^-$ corresponds (via the interpretation) to the $=$-class of $a$ in $\mcM\meq$
(so essentially the interpretation identifies every $a \in \mcM_P^-$ with itself). 

(iii) Suppose that
$0 < n < \omega$, $a_1, \ldots, a_n, b_1, \ldots, b_n \in M$, all $b_1, \ldots, b_n$ belong to the same
$E^P$-class and $b_i \neq a_i$ for all $i$. Assuming (without loss of generality) that $a_i \neq a_j$ if $i \neq j$ it follows that
$\{a_1, \ldots, a_n\}$ is an independent set (as $\mcM_P^-$ has SU-rank~1 and degenerate algebraic closure).
Moreover, $E^P(b_i, b_j)$ implies that 
\[
\tp_{(\mcM_P^-)\meq}(b_i / \acl_{(\mcM_P^-)\meq}(\es)) = 
\tp_{(\mcM_P^-)\meq}(b_j / \acl_{(\mcM_P^-)\meq}(\es))
\]
for all $i$ and $j$. Since $\mcM_P^-$ is $\omega$-categorical (and hence `small') it follows from
\cite[Corollary~5.3.5]{Kim_book} that 
$b_i$ and $b_j$ have the same Lascar strong type over $\es$, for all $i$ and $j$.
Hence the independence theorem for simple theories 
\cite[Theorem~3.2.8]{Kim_book} implies that there is $b$ such that 
$\tp_{\mcM_P^-}(a_i, b) = \tp_{\mcM_P^-}(a_i, b_i)$ for all $i = 1, \ldots, n$.

(iv) By Fact~\ref{basic facts about imaginaries}~(iii), $E^P$ is $\es$-definable in $\mcM_P^-$.
For a contradiction, suppose that $E^P$ is defined by a formula denoted $\varphi(x,y)$ and that no quantifier-free formula defines $E^P$.
Then we can find $a_1, a_2, b_1, b_2 \in M_P^-$ such that 
\[
\mcM_P^- \models \neg \varphi(a_1, a_2) \wedge \varphi(b_1, b_2), \ \ \text{ and } \ \ 
a_1a_2 \equiv_{\mcM_P^-}^\mrqf b_1b_2.
\]
From $a_1a_2 \equiv_{\mcM_P^-}^\mrqf b_1b_2$ and the definition of $\mcM_P^-$ it follows that
\[
\tp_\mcM^\mrqf(a_1, a_2 / A) = \tp_\mcM^\mrqf(b_1, b_2 / A).
\]
As $\mcM$ is homogeneous there is an automorphism of $\mcM$ which sends $a_1a_2$ to $b_1b_2$
and fixes $A$ pointwise.
By the definition of $\mcM_P^-$, the restriction of this automorphism to $M_P^-$ is an automorphism of $\mcM_P^-$,
so $a_1a_2 \equiv_{\mcM_P^-} b_1b_2$. Hence $\mcM_P^- \models \varphi(a_1, a_2) \wedge \varphi(b_1, b_2)$ or 
$\mcM_P^- \models \neg \varphi(a_1, a_2) \wedge \neg \varphi(b_1, b_2)$, but in either case it contradicts our assumption.
\hfill $\square$

\begin{defin}\label{definition of symmetric structure}{\rm
We say that a relational structure $\mcM$ is {\em symmetric} if for every relation symbol $R$ and
every $\bar{a} \in M$, $\mcM \models R(\bar{a})$ if and only if $\mcM \models R(\bar{a}')$ for every permutation $\bar{a}'$ of $\bar{a}$.
}\end{defin}

\begin{prop}\label{strong primitivity over a singleton}
Suppose that $\mcM$ is ternary, symmetric, 2-transitive, homogeneous, supersimple with SU-rank 1 and with degenerate algebraic closure.
If $b, b', a \in M$ and $\tp_\mcM(b / a) = \tp_\mcM(b' / a)$, then
$\tp_{\mcM\meq}(b / \acl_{\mcM\meq}(a)) = \tp_{\mcM\meq}(b' / \acl_{\mcM\meq}(a))$.
Consequently 
(by Lemma~\ref{degenerate acl, types over imaginaries and equivalence relations} and 
Theorem~\ref{weak isolation and nontrivial equivalence relations}) 
every constraint of cardinality 4 is weakly isolated.
\end{prop}

\noindent
{\bf Proof of Proposition~\ref{strong primitivity over a singleton}.}
Let $a \in M$, $A = \{a\}$, let $p$ be the unique nonalgebaric 1-type over $A$ and let $P = \{p\}$.
Also let $E^P$ be as in Proposition~\ref{M-p-minus is omega-categorical}, so $E^P$ has only infinite equivalence
classes and only finitely many equivalence classes.
It suffices to prove that $E^P$ has only one equivalence class.

Since $\mcM$ is 2-transitive it follows that if $a', b \in M$ are distinct and $q(x) = \tp_\mcM(b / a')$, then everything that we say about
$\mcM_P^-$ also holds for $\mcM_Q^-$ where $Q = \{q\}$.
Let 
\[
E^P_1(x, y), \ldots, E^P_k(x, y)
\] enumerate all 2-types (over $\es$) in $\mcM_P^-$ which imply that $E^P(x,y)$.
We will notationally identify each $E^P_i$ with the quantifier-free formula which isolates it.
Then $E^P(x, y)$ is equivalent, in $\mcM_P^-$, to $E^P_1(x, y) \vee \ldots \vee E^P_k(x, y)$.

Suppose that $E^P$ has at least two equivalence classes and let $X_1^P, X_2^P$ be two distinct $E^P$-classes.
The argument now splits into two cases, both of which will lead to a contradiction. 
Since the cases cover all possibilities it follows that $E^P$ has only one class.

\medskip

\noindent
{\em Case 1.}
Suppose that there are distinct $b, b', b'' \in M_P^-$ such that $b \in X_1^P$, $b', b'' \in X_2^P$ and $\mcM_Q^- \models E^Q(b, b'')$,
where $q$ is the unique nonalgebraic 1-type over $\{b'\}$ in $\mcM$ and $Q = \{q\}$.

\medskip

\noindent
From the symmetry of $\es$-definable ternary definable relations in $\mcM$ we get $b'b''a \equiv_\mcM b''ab'$.
Since $\mcM_P^- \models E^P(b', b'')$ it follows that $\mcM_Q^- \models E^Q(b'', a)$.
As $E^Q$ is transitive (being an equivalence relation) we get $\mcM_Q^- \models E^Q(b, a)$, so
$\mcM_Q^- \models E^Q_i(b, a)$ for some $i$.
Since $E^Q_i(x, y)$ isolates a type in $\mcM_Q^-$ it follows from the definition of $\mcM_Q^-$ that
$E^Q_i(x, y)$ determines the type of $xyb'$ in $\mcM$ (i.e. if $\mcM_Q^- \models E^Q_i(c, c') \wedge E^Q_i(d, d')$
then $cc'b' \equiv_\mcM dd'b'$).
By the symmetry of definable ternary relations in $\mcM$ we get
$bab' \equiv_\mcM bb'a$. 
Then $\mcM_P^- \models E^P_i(b, b')$ and hence  $\mcM_P^- \models E^P(b, b')$ which contradicts the choice of $b$ and $b'$.

\medskip

\noindent
{\em Case 2.}
Suppose that for all distinct $b, b', b'' \in M_P^-$, if $b \in X_1^P$ and $b', b'' \in X_2^P$,
then $\mcM_Q^- \models \neg E^Q(b, b'')$, where $q$ is the unique nonalgebraic 1-type over $\{b'\}$ in $\mcM$
and $Q = \{q\}$.

\medskip

\noindent
Given $b \in X_1^P$ let $r(x)$ be the unique nonalgebrac 1-type over $\{b\}$ in $\mcM$ and let $R = \{r\}$.
If $b', b'' \in X_2^P$ are distinct we have
$\mcM_Q^- \models \neg E^Q(b, b'')$ (where $Q$ depends on $b'$) and from the 
symmetry of $\es$-definable ternary relations in $\mcM$ we get 
$\mcM_R^- \models \neg E^R(b', b'')$.
Since we can fix (any) $b \in X_1^P$ and then vary $b', b''$ over all pairs of distinct elements in $X_2^P$ 
it follows, as $X_2^P$ is infinite, that $E^R$ has infinitely many equivalence classes.
But since $E^R$ has only infinite equivalence classes this contradicts that $\mcM$ has SU-rank 1.

Now we show that every constraint of cardinality 4 is weakly isolated.
Suppose for a contradiction that $\mcC$ is a constraint of cardinality 4 that is not weakly isolated.
Since $\mcM$ has (by assumption) degenerate algebraic closure it follows from
Theorem~\ref{weak isolation and nontrivial equivalence relations}
that there is $a \in M $, a nonalgebraic  $p \in S_1^\mcM(\{a\})$ and 
a nontrivial equivalence relation on $p(\mcM)$ which is $\{a\}$-definable and all of its equivalence classes are infinite.
By Lemma~\ref{degenerate acl, types over imaginaries and equivalence relations},
there are $b, b' \in M$ such that $\tp_\mcM(b / a) = \tp_\mcM(b' / a)$ and
\[
\tp_{\mcM\meq}(b / \acl_{\mcM\meq}(a)) \neq \tp_{\mcM\meq}(b' / \acl_{\mcM\meq}(a)).
\]
But this contradicts what we have proved.
\hfill $\square$

\begin{cor}\label{M-p-minus is a binary random structure}
Suppose that $\mcM$ is ternary, symmetric, 2-transitive, homogeneous, supersimple with SU-rank 1 and with
degenerate algebraic closure.
Let $a \in M$ and let $p(x)$ be the unique nonalgebraic 1-type of $\mcM$ over $\{a\}$.
Then $\mcM_{\{p\}}^-$ is a symmetric binary random structure.
\end{cor}

\noindent
{\bf Proof.}
By Proposition~\ref{M-p-minus is omega-categorical}, $\mcM_{\{p\}}^-$ is $\omega$-categorical, supersimple with SU-rank~1
and degenerate algebraic closure.
As $\mcM$ is symmetric it follows from the definition of $\mcM_{\{p\}}^-$ that $\mcM_{\{p\}}^-$ is symmetric.
Since $\mcM$ is 2-transitive it follows (by the definition of $\mcM_{\{p\}}^-$) that $\mcM_{\{p\}}^-$ has a unique 1-type over $\es$.
Now Proposition~\ref{strong primitivity over a singleton}
implies that $\mcM_{\{p\}}^-$ has a unique 1-type over $\acl_{(\mcM_{\{p\}}^-)\meq}(\es)$.
Hence $E^{\{p\}}$ (as in part~(iii) of Proposition~\ref{M-p-minus is omega-categorical}) has only one equivalence class.
It follows from Proposition~\ref{M-p-minus is omega-categorical}~(iii)
that if $0 < n < \omega$,
$a_1, \ldots, a_n, b_1, \ldots, b_n \in M_{\{p\}}^-$ and $b_i \neq a_i$ for all $i$, then there is $b \in M_{\{p\}}^-$ such that
$\tp_{\mcM_{\{p\}}^-}(a_i, b) = \tp_{\mcM_{\{p\}}^-}(a_i, b_i)$ for all $i = 1, \ldots, n$.
From this it follows that $\mcM_{\{p\}}^-$ has no constraint of cardinality greater than~2.
To show that $\mcM_{\{p\}}^-$ is a random structure 
(recall Definition~\ref{definition of age, permitted etc}~(v))
it remains to show that $\mcM_{\{p\}}^-$ is homogeneous.
But since $\mcM_{\{p\}}^-$ has a unique 1-type and has no constraint of cardinality greater than two it follows immediately that
its age has the disjoint amalgamation property and therefore $\mcM_{\{p\}}^-$ is homogeneous.
(Homogeneity can also be proved by a back-and-forth argument which builds up an automorphism that sends a given
tuple to any other given tuple that satisfies the same quantifier-free formulas.)
\hfill $\square$

\section{Examples}\label{Examples}

\noindent
Here we give examples of ternary 2-transitive homogeneous supersimple structures with SU-rank~1 and degenerate algebraic closure.
It remains open whether there is a structure with the same properties except that the algebraic closure is not degenerate (and hence, by
2-transitivity, nontrivial). 
The examples in Sections~\ref{Forbidding a complete 3-hypergraph} and~\ref{The parity 3-hypergraph}
are known.
Section~\ref{Uncountably many ternary homogeneous simple structures} shows
that even if the vocabulary $V$ contains only one relation symbol which is ternary
there are uncountably many homogeneous 2-transitive supersimple $V$-structures with SU-rank~1 and degenerate algebraic closure.
In Section~\ref{A ternary homogeneous structure which is definable in the generic tournament}
we give an example of a homogeneous ternary 2-transitive structure which is a reduct of the generic tournament.
I have not encountered the examples from Sections~\ref{Uncountably many ternary homogeneous simple structures}
and~\ref{A ternary homogeneous structure which is definable in the generic tournament}
in the literature or in the oral ``folklore''.\footnote{In \cite{AL} Akhtar and Lachlan show that there are
uncountably many homogeneous 3-hypergraphs, but it seems like their examples are not simple although I have not checked this.}

All examples in Sections~\ref{Forbidding a complete 3-hypergraph}
--~\ref{A ternary homogeneous structure which is definable in the generic tournament} 
have only weakly isolated constraints.
Also, the same examples split into two categories: the examples with the free amalgamation property;
every other example (i.e. those in Sections~\ref{The parity 3-hypergraph}
and~\ref{A ternary homogeneous structure which is definable in the generic tournament})
is a reduct of a binary random structure.

If, for some vocabulary $V$, $\mbC$ is a class of finite $V$-structures, then $\mbF(\mbC)$ denotes the class of all
finite $V$-structures $\mcA$ such that no member of $\mbC$ embeds into $\mcA$. If $\mbC = \{\mcC_1, \ldots, \mcC_k\}$, then
we may write $\mbF(\mcC_1, \ldots, \mcC_k)$ instead of $\mbF(\mbC)$.
We will consider examples in which we consider only structures in which some relation symbol is always interpreted as a symmetric relation
and in these cases we do not explicitly mention the constraints which express this.

\subsection{Forbidding 3-irreducible structures}\label{Forbidding a complete 3-hypergraph}

The notion of  `indecomposable structure' used by Henson in \cite{Hen72} has been generalized by Conant to the notion of 
`$k$-irreducible structure' in \cite{Con}.
We say that a structure $\mcA$ (in any relational language) is {\em $k$-irreducible} if for any choice of $k$ elements $a_1, \ldots, a_k \in A$
there is a relationship $\bar{a}$ (of $\mcA$) such that $a_1, \ldots, a_k \in \rng(\bar{a})$.
We will use the following two results.

\begin{fact}\label{on indecomposable constraints} {\rm (Henson \cite[Theorem 1.2]{Hen72})}
For any finite relational vocabulary $V$, if $\mbC$ is
a class of finite $V$-structures and every member of $\mbC$ is 2-irreducible, then $\mbF(\mbC)$ has the free amalgamation property.
\end{fact}

\noindent
We recall that if $\mbK$ has the disjoint amalgamation property (and the hereditary property) then the
Fra\"{i}ss\'{e} limit of $\mbK$ has degenerate algebraic closure.
Part~(a) of the following result shows how to construct many ternary homogeneous supersimple structures with SU-rank~1
and degenerate algebraic closure. It will be used more systematically in 
Section~\ref{Uncountably many ternary homogeneous simple structures}.
If $\mcA$ and $\mcB$ are relational structures with the same vocabulary
and $f : A \to B$ we call $f$ a {\em homomorphism} from $\mcA$ to $\mcB$ if for every relation symbol $R$ and every $\bar{a} \in A$,
$\mcA \models R(\bar{a})$ implies $\mcB \models R(f(\bar{a}))$.

\begin{fact}\label{Conants result on 3-irreducible constraints} {\rm (Conant \cite[Theorem 7.22]{Con})}
Let $V$ be a finite relational vocabulary and suppose that $\mbC$ is a set of finite $V$-structures such for any two different
$\mcA, \mcB \in \mbC$ there is no injective homomorphism from $\mcA$ into $\mcB$.
Let $\mbK$ be the class of all finite $V$-structures $\mcA$ such that for every $\mcC \in \mbC$ there does not exist an
injective homomorphism from $\mcC$ into $\mcA$.
\begin{itemize}
\item[(a)] If every structure in $\mbC$ is 3-irreducible then the Fra\"{i}ss\'{e} limit of $\mbK$ is supersimple with 
SU-rank~1 and degenerate algebraic closure.
\item[(b)] If $\mbK$ has the free amalgamation property, then its Fra\"{i}ss\'{e} limit is simple if and only if every member of 
$\mbC$ is 3-irreducible.
\end{itemize}
\end{fact}

\begin{rem}\label{remark on conants result}{\rm
(a) Part~(a) of Fact~\ref{Conants result on 3-irreducible constraints} holds (by inspection of its proof)
if  `$\mbK$' is replaced by `$\mbF(\mbC)$' and the assumption on $\mbC$ is changed as follows:
suppose that for any two different $\mcA, \mcB \in \mbC$, there is no embedding from $\mcA$ into $\mcB$.
(This observation is due to Gabriel Conant.) \\
(b) If we want to consider a class of finite structures in which some relation symbol, say $R$, is always interpreted as
a symmetric relation, then, in order to use Fact~\ref{Conants result on 3-irreducible constraints}, we must
{\em not} add constraints to $\mbC$ which express this. Instead we just ignore all structures in which $R$ is not interpreted as
a symmetric relation.
To see why, consider this example:
The structure $\mcA = (\{a, b, c\}, R^\mcA)$, where $R^\mcA = \{(a, b, c)\}$, does not belong to the class of all
finite 3-hypergraphs (viewed as $\{R\}$-structures) because $R^\mcA$ is not a symmetric relation, 
but there is an injective homomorphism from $\mcA$ into any
3-hypergraph with at least one edge. 
Consequently, the class of all finite 3-hypergraphs which do not have a complete 3-hypergraph on 4 vertices as a subgraph
cannot be described in the framework of Fact~\ref{Conants result on 3-irreducible constraints},
{\em unless} we simply ignore all $\{R\}$-structures in which $R$ is not interpreted as a symmetric relation.

(c) Note that if $V$ is ternary and $\mbK$ is as in Fact~\ref{Conants result on 3-irreducible constraints}
and every member of $\mbC$ is 3-irreducible,
then every constraint of $\mbK$ is weakly isolated.
}\end{rem}

\noindent
A particularly well known example that can be obtained from 
Fact~\ref{Conants result on 3-irreducible constraints}~(a) with the interpretation of
part~(b) of the above remark is the generic tetrahedron-free 3-hypergraph, that is, the 
Fra\"{i}ss\'{e} limit of the class of all finite 3-hypergraphs into which $\mcK_4$ cannot be embedded, where
$\mcK_4$ denotes the complete 3-hypergraph on 4 vertices.
As a contrast consider the following example.
Let $\mcK_4^-$ be the result of removing one hyperedge from $\mcK_4$.
Since $\mcK_4^-$ and $\mcK_4$ are 2-irreducible, $\mbF(\mcK_4^-, \mcK_4)$ has the free amalgamation property by 
Fact~\ref{on indecomposable constraints}.
Since $\mcK_4^-$ is not 3-irreducible it follows from part~(b) of Fact~\ref{Conants result on 3-irreducible constraints} that the 
Fra\"{i}ss\'{e} limit of $\mbF(\mcK_4^-, \mcK_4)$ is not simple.
But it is superrosy with thorn rank~1, which can be concluded from results in \cite{Con}.

\subsection{The ``parity 3-hypergraph''}\label{The parity 3-hypergraph}

In \cite{AL} Akhtar and Lachlan study infinite homogeneous 3-hypergraphs.
One of the examples that they consider is the following one, which they attribute to Cherlin and Macpherson. 
Let $V = \{R\}$ and $V_E = \{E\}$ where $E$ is a binary and $R$ a ternary relation symbol.
Let $\mcG$ be the Rado graph, with vertex set $G$, viewed as a $V_E$-structure.
Furthermore, let $\mcH$ be the 3-hypergraph, viewed as a $V$-structure with the same vertex set as $\mcG$ (i.e. $H = G$) where, for
all distinct vertices $a, b$ and $c$, $\mcH \models R(a, b, c)$ if and only if the number of edges (of $\mcG$) between elements in the set
$\{a, b, c\}$ is odd. 
As explained in \cite{AL}, $\mcH$ is homogeneous. Since $\mcH$ is a reduct of $\mcG$ and
$\mcG$ is supersimple with SU-rank 1 and with degenerate algebraic closure, 
the same is true for $\mcH$ (by Fact~\ref{facts about interpretability}).

The structure $\mcH$ can also be characterized in the following way (and we refer to \cite{AL} for more explanations
of the nontrivial claims that follow).
Let $\mcC_1$ and $\mcC_3$ be 3-hypergraphs with exactly four vertices and such that 
$\mcC_1$ has exactly 1 hyperedge and $\mcC_3$ has exactly 3 hyperedges.
Then $\mbF(\mcC_1, \mcC_3)$ has the disjoint amalgamation property (but not the free amalgamation property)
and $\mcH$ is the Fra\"{i}ss\'{e} limit of $\mbF(\mcC_1, \mcC_3)$.
It is immediate that both $\mcC_1$ and $\mcC_3$ are isolated with respect to $\mbF(\mcC_1, \mcC_3)$.

\subsection{Uncountably many ternary homogeneous simple structures}\label{Uncountably many ternary homogeneous simple structures}

Let $V = \{R\}$ where $R$ is a ternary relation symbol.
With the help of Fact~\ref{Conants result on 3-irreducible constraints}
we modify, for our purposes, the construction which Henson used to obtain uncountably many homogeneous directed graphs \cite{Hen72}.
For each $2 < n < \omega$, let $\mcH_n$ be the $V$-structure with universe $\{0, \ldots, n\}$ such that
\begin{itemize}
\item if $\mcH_n \models R(a, b, c)$ then $a$, $b$ and $c$ are distinct, and
\item for all distinct $a, b, c \in \{0, \ldots, n\}$, 
\begin{itemize}
\item[] $\mcH_n \models \neg R(a, b, c)$ if and only if $a = 0$, $b > 0$ and either 
$b < n$ and $c = b+1$, or $b = n$ and $c = 1$.
\end{itemize}
\end{itemize}

\begin{lem}\label{H-n cannot be embedded into H-m}
If $n \neq m$ then there is no embedding from $\mcH_n$ to $\mcH_m$.
\end{lem}

\noindent
{\bf Proof.}
Let $n < m$.
Suppose for a contradiction that $f : \mcH_n \to \mcH_m$ is an embedding.
First note that in both $\mcH_n$ and $\mcH_m$ all elements but 0 satisfy the following formula:
\[
\forall y, z \big( x \neq y \wedge y \neq z \wedge x \neq z \ \rightarrow \ R(x, y, z)\big).
\]
Hence $f(0) = 0$.
Observe that the map $g : \{0, \ldots, n\} \to \{0, \ldots, n\}$ given by $g(0) = 0$,
$g(k) = k+1$ if $0 < k < n$ and $g(n) = 1$ is an automorphism of $\mcH_n$.
Therefore we may,  without loss of generality, assume that $f(1) = 1$.
By the definition of $\mcH_n$, for each $k \in \{1, \ldots, n\}$ there is a unique $l \in \{1, \ldots, n\}$
such that $\mcH_n \models \neg R(0, k, l)$ (and we must have either $l = k+1$ or $l = 1$ and $k = n$).
The same is true if $n$ is replaced by $m$.
It follows that we must have $f(k) = k$ for all $k \in \{1, \ldots, n\}$ and consequently $n = m$ which contradicts our assumption.
\hfill $\square$
\\

\noindent
Let $\mbS = \{\mcH_n : 2 < n < \omega\}$ and $\mbT \subseteq \mbS$.
Clearly, for every $2 < n <\omega$, $\mcH_n$ is 3-irreducible and hence 2-irreducible.
So by Fact~\ref{on indecomposable constraints}, $\mbF(\mbT)$ has the free amalgamation property and
(since it obviously has the hereditary property)
we can let $\mcM_\mbT$ denote the Fra\"{i}ss\'{e} limit of $\mbF(\mbT)$.
From Lemma~\ref{H-n cannot be embedded into H-m},
Fact~\ref{Conants result on 3-irreducible constraints}~(a) and
Remark~\ref{remark on conants result}~(a)
it follows that
$\mcM_\mbT$ is supersimple with SU-rank~1 and degenerate algebraic closure.
From Lemma~\ref{H-n cannot be embedded into H-m} again,
it follows that if $T$ and $T'$ are different subsets of $\mbS$, then  $\mcM_\mbT \not\cong \mcM_{\mbT'}$.
As $\mbT \subseteq \mbS$ can be chosen in $2^\omega$ ways, we get $2^\omega$ nonisomorphic (countable) structures $\mcM_\mbT$.
It follows directly from the definitions that for every $\mbT \subseteq \mbS$ and every $\mcH_n \in \mbC_\mbT$,
$\mcH_n$ is isolated with respect to $\mcM_\mbT$.

\subsection{A ternary homogeneous structure which is a reduct of the generic tournament}
\label{A ternary homogeneous structure which is definable in the generic tournament}

Let $E$ be a binary relation symbol and let $\mbD$ be the class of all finite tournaments, viewed as $\{E\}$-structures.
It is easy to see that $\mbD$ has the hereditary property and the (disjoint, but not free) amalgamation property.
Let $\mcG$ be the Fra\"{i}ss\'{e} limit of $\mbD$.

\begin{defin}\label{definition of n-reversal equivalence}{\rm
Let $1 < n < \omega$ and $(a_1, \ldots, a_n), (b_1, \ldots, b_n)  \in G^n$. 
We define 
\[
(a_1, \ldots, a_n) \approx_n (b_1, \ldots, b_n)
\] 
if and and only if the following two conditions hold:
\begin{itemize}
\item[(a)] For all $1 \leq i, j \leq n$, $a_i = a_j$ $\Longleftrightarrow$ $b_i = b_j$.

\item[(b)] Either
\begin{itemize}
\item[] for all $1 \leq i, j \leq n$, $\mcG \models E(a_i, a_j) \leftrightarrow E(b_i, b_j)$, 
\item[] or, for all $1 \leq i, j \leq n$, $\mcG \models E(a_i, a_j) \leftrightarrow E(b_j, b_i)$.
\end{itemize}
\end{itemize}
}\end{defin}

\noindent
Observe that `$\approx_n$' is an equivalence relation which is $\es$-definable in $\mcG$.
Also note that for all pairs of distinct elements $(a_1, a_2), (b_1, b_2) \in G^2$ we have $(a_1, a_2) \approx_2 (b_1, b_2)$.

\begin{lem}\label{n-reversal equivalence and 3-reversal equivalence}
Let $n \geq 3$ and $(a_1, \ldots, a_n), (b_1, \ldots, b_n) \in G^n$.
Then $(a_1, \ldots, a_n) \approx_n (b_1, \ldots, b_n)$ if and only if
for all distinct $i, j, k \in \{1, \ldots, n\}$,
 $(a_i, a_j, a_k) \approx_3 (b_i, b_j, b_k)$.
\end{lem}

\noindent
{\bf Proof.}
Only one of the implications is nontrivial.
Let $(a_1, \ldots, a_n), (b_1, \ldots, b_n) \in G^n$ where $n \geq 3$ and suppose that 
$(a_i, a_j, a_k) \approx_3 (b_i, b_j, b_k)$ for all distinct $i, j, k \in \{1, \ldots, n\}$.
By induction on $n$ we prove that $(a_1, \ldots, a_n) \approx_n (b_1, \ldots, b_n)$.
The base case is $n = 3$ and then there is nothing to prove.

So suppose that $n > 3$.
By the induction hypothesis, for every $i \in \{1, \ldots, n\}$,
\begin{equation}\label{the induction hypothesis to approx-n-1}
(a_1, \ldots, a_{i-1}, a_{i+1}, \ldots, a_n) \approx_{n-1} (b_1, \ldots, b_{i-1}, b_{i+1}, \ldots, b_n).
\end{equation}
In particular we have $(a_1, \ldots, a_{n-1}) \approx_{n-1} (b_1, \ldots, b_{n-1})$, so either for all
distinct $i, j \in \{1, \ldots, n-1\}$,
$\mcG \models E(a_i, a_j) \leftrightarrow E(b_i, b_j)$, 
or for all distinct $i, j \in \{1, \ldots, n-1\}$,
$\mcG \models E(a_i, a_j) \leftrightarrow E(b_j, b_i)$.
Accordingly we divide the argument into two (similar) cases.

\medskip

\noindent
{\em Case 1: For all distinct $i, j \in \{1, \ldots, n-1\}$, $\mcG \models E(a_i, a_j) \leftrightarrow E(b_i, b_j)$}

\medskip

\noindent
Let $i, j \in \{1, \ldots, n-1\}$ be distinct.
Since $\mcG$ is a tournament we can, without loss of generality, assume that $\mcG \models E(a_i, a_j) \wedge E(b_i, b_j)$.
First suppose that $\mcG \models E(a_i, a_n) \wedge E(a_j, a_n)$.
Since $(a_i, a_j, a_n) \approx_3 (b_i, b_j, b_n)$ it follows that $\mcG \models E(b_i, b_n) \wedge E(b_j, b_n)$.
If instead $\mcG \models E(a_n, a_i) \wedge E(a_j, a_n)$, then, by the same argument, we see that
$\mcG \models E(b_n, b_i) \wedge E(b_j, b_n)$.
We  argue in the same way in the remaining two cases.
As the arguments works for any distinct $i, j \in \{1, \ldots, n-1\}$ it follows that for all $1 \leq i \leq n-1$, 
$\mcG \models E(a_i, a_n) \leftrightarrow E(b_i, b_n)$ and therefore
$(a_1, \ldots, a_n) \approx_n (b_1, \ldots, b_n)$.

\medskip

\noindent
{\em Case 2: For all distinct $i, j \in \{1, \ldots, n-1\}$, $\mcG \models E(a_i, a_j) \leftrightarrow E(b_j, b_i)$}

\medskip

\noindent
Let $i, j \in \{1, \ldots, n-1\}$ be distinct.
Since $\mcG$ is a tournament we can, without loss of generality, assume that $\mcG \models E(a_i, a_j) \wedge E(b_j, b_i)$.
The rest of the argument in Case~2 is an obvious modification of the argument in Case~1, which is left for the reader.
\hfill $\square$
\\

\noindent
It is straightforward to verify that $\approx_3$ has exactly four equivalence classes, say $X_1, \ldots, X_4$, on triples of {\em distinct} elements.
Let $V = \{R_1, \ldots, R_4\}$ where $R_1, \ldots, R_4$ are ternary relation symbols.
Let  $\mcM$ be the $V$-structure with the same universe as $\mcG$ (so $M = G$)  such that for each  $i = 1, \ldots, 4$,
$R_i^\mcM = X_i$.

\begin{lem}\label{the reversal structure is homogeneous}
$\mcM$ is homogeneous and supersimple with SU-rank 1 and with degenerate $\acl_\mcM$.
\end{lem}

\noindent
{\bf Proof.}
From the definitions of $\approx_3$ and $R_i^\mcM$, each $R_i^\mcM$ is $\es$-definable in $\mcG$.
Since $\mcG$ is supersimple with SU-rank 1 and $\acl_\mcG$ is degenerate
it follows (using Fact~\ref{facts about interpretability}) that the same is true for $\mcM$ and $\acl_\mcM$.
So it remains to show that $\mcM$ is homogeneous.

Let $0 < n < \omega$, $a_1, \ldots, a_n, a_{n+1}, b_1, \ldots, b_n \in M = G$ and suppose that 
\begin{equation}\label{bar-a and bar-b have the same quantifier-free type}
\tp_\mcM^\mrqf(a_1, \ldots, a_n) = \tp_\mcM^\mrqf(b_1, \ldots, b_n).
\end{equation}
We need to find $b_{n+1} \in M$ such that 
\[
\tp_\mcM^\mrqf(a_1, \ldots, a_n, a_{n+1}) = \tp_\mcM^\mrqf(b_1, \ldots, b_n, b_{n+1}).
\]
By the definition of $\mcM$, this is accomplished if we find $b_{n+1} \in G = M$ such that 
$(a_1, \ldots, a_n, a_{n+1}) \approx_n (b_1, \ldots, b_n, b_{n+1})$.

We consider only the case when $n > 2$ since the cases $n = 1, 2$ are similar and simpler.
From the definition of $\mcM$ and~(\ref{bar-a and bar-b have the same quantifier-free type}) 
it follows that for all distinct $i, j, k \in \{1, \ldots, n\}$,
$(a_i, a_j, a_k) \approx_3 (b_i, b_j, b_k)$.
By Lemma~\ref{n-reversal equivalence and 3-reversal equivalence}
we get 
\begin{equation}\label{a-bar and b-bar are n-reversal equivalent}
(a_1, \ldots, a_n) \approx_n (b_1, \ldots, b_n).
\end{equation}
Without loss of generality we can assume that $a_{n+1} \notin \{a_1, \ldots, a_n\}$.
In order to find $b_{n+1} \in G$ such that 
$(a_1, \ldots, a_n, a_{n+1}) \approx_n (b_1, \ldots, b_n, b_{n+1})$,
it suffices, by the definition of $\approx_n$ and~(\ref{a-bar and b-bar are n-reversal equivalent}), to find $b_{n+1} \in G$ such that 

\begin{itemize}
\item[] if there are distinct $i, j \in \{1, \ldots, n\}$ such that 
$\mcG \models E(a_i, a_j) \wedge E(b_i, b_j)$ (which, via~(\ref{a-bar and b-bar are n-reversal equivalent}), implies that for all 
$i, j \in \{1, \ldots, n\}$, $\mcG \models E(a_i, a_j) \leftrightarrow E(b_i, b_j)$), then 
\[
\mcG \models E(a_i, a_{n+1}) \leftrightarrow E(b_i, b_{n+1}) \  \text{ for all } i \in \{1, \ldots, n\},
\]
and otherwise
\[
\mcG \models E(a_i, a_{n+1}) \leftrightarrow E(b_{n+1}, b_i) \ \text{ for all } i \in \{1, \ldots, n\}.
\]
\end{itemize}
Since $\mcG$ is the generic tournament (i.e. the Fra\"{i}ss\'{e} limit of $\mbD$) it is possible to find such $b_{n+1} \in G$.
\hfill $\square$

\begin{lem}\label{the reversal structure has only finitely many constraints}
(i) Every constraint of $\mcM$ has at most 4 elements.\\
(ii) Every constraint of $\mcM$ is weakly isolated.\\
(iii) The age of $\mcM$ does not have the free amalgamation property.
\end{lem}

\noindent
{\bf Proof.}
(i) Suppose, towards a contradiction, that $\mcC$ is a constraint of $\mcM$ with at least 5 elements.
Let $c \in C$ and let $\mcA = \mcC \uhrc C \setminus \{c\}$, so $\mcA$ is permitted with respect to $\mcM$.
Then we may, without loss of generality, assume that $\mcA \subseteq \mcM$.
As $\mcC$ is a constraint with at least 5 elements it follows that for every triple $\bar{a} = (a_1, a_2, a_3) \in A^3$ 
(where  $A = C \setminus \{c\}$) there
is $c_{\bar{a}} \in M$ such that $\tp_\mcC^\mrqf(\bar{a}, c) = \tp_\mcM^\mrqf(\bar{a}, c_{\bar{a}})$.

Suppose that
$\bar{a} = (a_1, a_2, a_3), \bar{a}' = (a_1, a_2, a'_3) \in A^3$ are triples of distinct elements.
From the choice of $c_{\bar{a}}$ and $c_{\bar{a}'}$ it follows that
$\tp_\mcM^\mrqf(a_1, a_2, c_{\bar{a}}) = \tp_\mcM^\mrqf(a_1, a_2, c_{\bar{a}'})$, so
$(a_1, a_2, c_{\bar{a}}) \approx_3 (a_1, a_2, c_{\bar{a}'})$ and hence
$\tp_\mcG^\mrqf(a_1, a_2, c_{\bar{a}}) = \tp_\mcG^\mrqf(a_1, a_2, c_{\bar{a}'})$.

Now suppose that $\bar{a} = (a_1, a_2, a_3), \bar{a}' = (a_1, a'_2, a'_3) \in A^3$ are triples of distinct elements
where $\{a_2, a_3\} \cap \{a'_2, a'_3\} = \es$.
Let $\bar{a}'' = (a_1, a_2, a'_2) $. From what we just proved it follows that
\[
\tp_\mcG^\mrqf(a_1, c_{\bar{a}}) = \tp_\mcG^\mrqf(a_1, c_{\bar{a}''}) = \tp_\mcG^\mrqf(a_1, c_{\bar{a}'}).
\]
Thus we have proved that if $\bar{a} = (a_1, a_2, a_3), \bar{a}' = (a_1, a'_2, a'_3) \in A^3$ are triples of distinct elements,
then $\tp_\mcG^\mrqf(a_1, c_{\bar{a}}) = \tp_\mcG^\mrqf(a_1, c_{\bar{a}'})$.
Since $\mcG$ is the generic tournament it follows that
\[
\big\{ \varphi(x, \bar{a}) : \varphi(x, \bar{a}) \in \tp_\mcG^\mrqf(c_{\bar{a}} / \bar{a}) \text{ and } \bar{a} \in A^3 \big\}
\]
is consistent. So there is $c' \in M$ such that for all $\bar{a} \in A^3$,
$\tp_\mcG^\mrqf(c', \bar{a}) = \tp_\mcG^\mrqf(c_{\bar{a}}, \bar{a})$.
Then $\tp_\mcM^\mrqf(c', \bar{a}) = \tp_\mcM^\mrqf(c_{\bar{a}}, \bar{a})$ for every $\bar{a} \in A^3$
and consequently
$\tp_\mcM^\mrqf(c' / A) = \tp_\mcC^\mrqf(c / A)$ which contradicts that $\mcC$ is a constraint.

(ii) As $\mcG$ is the generic tournament, it is easy to see that for every finite $A \subseteq G$ there is no nontrivial
$A$-definable (in $\mcG$) equivalence relation on $G \setminus A$.
Since $\mcM$ is a reduct of $\mcG$ it follows that the same holds for $\mcM$.
Hence the conclusion follows from Corollary~\ref{no nontrivial equivalence relations imply only weakly isolated constraints}.

Part~(iii) follows from the fact that for every triple $(a_1, a_2, a_3)$ of distinct elements from $M$ we have,
by the definition of $\mcM$,
$\mcM \models \bigvee_{i = 1}^4 R_i(a_1, a_2, a_3)$. 
\hfill $\square$

\subsection{Examples of higher SU-rank}\label{Examples of higher SU-rank}

All examples above have SU-rank 1. 
It is not hard, however, to construct examples of higher SU-rank, for example 2, by using a binary relation symbol interpreted as
an equivalence relation with infinitely many infinite equivalence classes. If one likes, on each class 
(which in itself is a structure of rank 1) one can add some more ``exotic'' structure by (for example) 
choosing one of the structures earlier in this section, call it $\mcM$, and letting each class, as a structure in itself be
isomorphic to $\mcM$. Obviously this kind of example is not primitive, i.e. it has a nontrivial $\es$-definable equivalence relation
on its universe.

A different kind of example of SU-rank 2 which is primitive is Example~3.3.2 in~\cite{Mac11}, which is also discussed in 
\cite[Example~2.7]{KopPrimitive}. This example is not 2-transitive.
In fact, it follows from Observation~\ref{observation about 2-transitivity and triviality}
that if $\mcM$ is a homogeneous simple structure with trivial dependence and higher SU-rank than one, then it
must have some binary relation symbol (it is easy to see that unary relation symbols would not suffice to raise the rank).

\section{Problems}\label{Problems}

\noindent
As a number of questions are left unanswered, also taking the conclusions of the later article
\cite{Kop17b} into account, we conclude with a collection of problems.

\begin{enumerate}

\item Is every ternary homogeneous  simple structure supersimple (with finite SU-rank)?

\item For $k \geq 4$, is every $k$-ary homogeneous finitely constrained simple structure supersimple (with finite SU-rank)?

\item Is there $k \geq 4$ and a $k$-ary homogeneous supersimple structure $\mcM$ with SU-rank~1 and 
nontrivial algebraic closure?

\item Suppose that $\mcM$ is ternary, 2-transitive, homogeneous, supersimple with SU-rank~1 and degenerate algebraic closure.
\begin{enumerate}

\item Can $\mcM$ have a constraint which is not weakly isolated?

\item Can $\mcM$ have a finite subset $A \subset M$, a nonalgebraic type $p \in S_1^\mcM(A)$ and a nontrivial
$A$-definable equivalence relation on $p(\mcM)$?
(A negative answer implies, by Theorem~\ref{weak isolation and nontrivial equivalence relations}, 
a negative answer to the previous question.)

\item If there is a finite $A \subset M$, a nonalgebraic type $p \in S_1^\mcM(A)$ and a nontrivial $A$-definable
equivalence relation on $p(\mcM)$, does it follow that $\mcM$ has a constraint which is not weakly isolated?

\item If the age of $\mcM$ does {\em not} have the free amalgamation property, must $\mcM$ be a reduct of
(or more generaly, interpretable in) a binary random structure?
\end{enumerate}
\end{enumerate}

\noindent
{\bf Acknowledgement.} I thank the anonymous referee for a detailed examination of the article,
including finding some minor mistakes (now corrected) and
giving suggestions that improved the clarity of the arguments.

\end{document}